\documentclass[12pt]{article}
\input epsf
\usepackage{amsfonts}
\usepackage{amscd}
\usepackage{amsmath,amssymb}
\usepackage{amstext}
\usepackage{amscd,graphicx}
\usepackage[all]{xy}

\newtheorem{lem}{LEMMA}[section]
\newtheorem{theo}[lem]{THEOREM}
\newtheorem{cla}[lem]{Claim}
\newtheorem{coro}[lem]{COROLLARY}
\newtheorem{prop}[lem]{PROPOSITION}

\newtheorem{definition}[lem]{Definition}

\newtheorem{rem}[lem]{Remark}
\newtheorem{ex}[lem]{Example}

\renewcommand{\descriptionlabel}[1]%
       {\hspace{\labelsep}\textsf{#1}}

\begin{document}

\title{A wild knot $\mathbb{S}^{2}\hookrightarrow\mathbb{S}^{4}$ as
  limit set of a Kleinian Group: Indra's pearls in four dimensions.}
\author{Gabriela Hinojosa \thanks{This work was partially supported by
CONACyT (Mexico), grant G36357-E, and DGEP-UNAM (M\'exico).}\\
Facultad de Ciencis, UAEM\\
62210, Cuernavaca, M\'exico.\\
gabriela@matcuer.unam.mx}
\date{July 31, 2003}
\maketitle
\begin{abstract}
The purpose of this paper is to construct an example of a 2-knot
wildly embedded in $\mathbb{S}^{4}$ as the limit set of a Kleinian
group. We find that this type of wild 2-knots has very interesting
topological properties. 
 \end{abstract}

\section{Introduction}

One geometrical method used to obtain remarkable fractal sets of
extreme beauty and complexity having the property of being
self-similar (i.e. conformally equal to itself at infinite small scales), is by
considering the limit sets of Schottky groups, consisting on finitely
generated groups of reflections on codimension one round spheres. As a
testimony of such a beauty and complexity, one can consult the wonderful
book {\it Indra's Pearl: The vision of Felix Klein} written by
D. Mumdord, C. Series and D. Wright \cite{mumford}.\\

The purpose of the present paper is to construct, in the spirit of
Indra's pearls book, an example of a wildly embedded 2-sphere in
$\mathbb{S}^{4}$ (i.e. a wild 2-knot in $\mathbb{S}^{4}$) obtained as
limit set of a Kleinian group.\\

In section 2, we present the preliminary definitions and results in
knot theory and Kleinian groups that we will use in this paper. In
section 3, we describe the geometric ideas involved to construct a
wild 2-knot, and we give an explicit example of such a group. In
section 4, we prove that the limit set obtained in section 3, is a wild
2-knot in $\mathbb{S}^{4}$. In sections 5, 6 and 7 we give very
interesting topological properties in the case where the original arc (see section
2) fibers over the circle. We show that the wild 2-knot also fibers over the
circle and we determine its monodromy. In section 8, we
lift the action of this Kleinian group to the
twistor space of $\mathbb{S}^{4}$, obtaining a dynamically defined 
$\mathbb{S}^{2}\times \mathbb{S}^{2}$ wildly embedded in the twistor space 
$P^{3}_{\mathbb{C}}$.\\

I want to thank Prof. Alberto Verjovsky for all the discussions
and very important suggestions. I also want to thank Prof. Aubin
Arroyo for drawing the beautiful Figures 10 and 11.\\ 

\section{Preliminaries}

In 1925 Emil Artin described two methods for constructing knotted spheres of
dimension two in $\mathbb{S}^{4}$ from knots in $\mathbb{S}^{3}$. The first
of them is called {\it suspension}. Roughly speaking, this method consists
of taking the suspension of $(K,\mathbb{S}^{3})$, where $K\subset\mathbb{S}^{3}$ is a 
tame knot, to obtain a 2-knot $\Sigma K$ in $\mathbb{S}^{4}$. By construction, we have that 
the fundamental group of $\Sigma K$ is easily computable.\\

The second method is called {\it spinning} and  uses the rotation process.
A way to visualize it is the following. We can consider $\mathbb{S}^{2}$ as an
$\mathbb{S}^{1}$-family of half-equators (meridians=$\mathbb{D}^{1}$) such that
the respective points of their boundaries are identified to obtain the poles. Then the formula 
$Spin(\mathbb{D}^{1})=\mathbb{S}^{2}$ means to send 
homeomorphically the unit interval $\mathbb{D}^{1}$ to a meridian of
$\mathbb{S}^{2}$ such that $\partial\mathbb{S}^{1}=\{0,1\}$ is mapped to the poles
and, multiply the interior of $\mathbb{D}^{1}$ by $\mathbb{S}^{1}$. In other words, one 
spins the meridian with respect to the poles to obtain $\mathbb{S}^{2}$. 
Similarly, consider $\mathbb{S}^{n+1}$ as an $\mathbb{S}^{1}$-family of half-equators 
($\mathbb{D}^{n}$) where boundaries are respectively identified, hence 
$Spin(\mathbb{D}^{n})=\mathbb{S}^{n+1}$ means to send homeomorphically 
$\mathbb{D}^{n}$ to a meridian of $\mathbb{S}^{n+1}$ and keeping $\partial\mathbb{D}^{n}$
fixed, multiply the interior of $\mathbb{D}^{n}$ by $\mathbb{S}^{1}$. In 
particular $Spin(\mathbb{D}^{3})=\mathbb{S}^{4}$.\\

It is this second method that we will use to construct a 2-sphere wildly embedded
in $\mathbb{S}^{4}$, so we will give a more detailed description of it.\\

Consider in $\mathbb{R}^{4}$ the half-space
$$
\mathbb{R}^{3}_{+}=\{(x_{1},x_{2},x_{3},0):x_{3}\geq 0\}
$$
whose boundary is the plane
$$
\mathbb{R}^{2}=\{(x_{1},x_{2},0,0)\}.
$$
We can spin each point $x=(x_{1},x_{2},x_{3},0)$ of $\mathbb{R}^{3}_{+}$
with respect to $\mathbb{R}^{2}$ according to the formula
$$
R_{\theta}(x)=(x_{1},x_{2},x_{3}\cos\theta,x_{3}\sin\theta).
$$
We define $Spin(X)$ of a set $X\subset\mathbb{R}^{3}_{+}$, as
$$
Spin(X)=\{R_{\theta}(x):x\in X, 0\leq\theta\leq 2\pi\}.
$$
To obtain a knot in $\mathbb{R}^{4}$, we choose a tame arc $A$ in 
$\mathbb{R}^{3}_{+}$ with its end-points in $\mathbb{R}^{2}$ and its interior in 
$\mathbb{R}^{3}_{+}\setminus\mathbb{R}^{2}$. Then
$Spin(A)$ is a 2-sphere in $\mathbb{R}^{4}$ called a {\it spun knot}.\\

\begin{figure}[tbh]
\centerline{\epsfxsize=2in \epsfbox{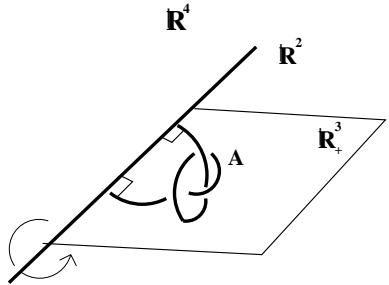}}
\caption{\sl Spun knot.}
\label{F1}
\end{figure}

We can think of $A$ as the image of an embedding 
$A:I\rightarrow\mathbb{R}^{3}_{+}$ with $A(0)\neq A(1)\in\mathbb{R}^{2}$, and which will
be denoted by the same letter. Then we will say that an arc $A\subset\mathbb{R}^{3}_{+}$
is a {\it spinnable arc} if it is smooth in every point with contact of infinite order with
respect to the normal on its end-points.\\
 
It can be proved that the fundamental group of $Spin(A)$ is isomorphic to 
$\Pi_{1}(\mathbb{R}^{3}_{+}-A)$ and by the Seifert-Van Kampen
Theorem's, this is isomorphic to the fundamental group
of $A\cup L$ in $\mathbb{R}^{3}$, where $L\subset\mathbb{R}^{3}$ is an unknotted
segment joining the end-points of $A$ (see \cite{rolfsen}, \cite{zeeman1}).\\

Our goal is to obtain a wild 2-sphere as the limit set of a conformal 
Kleinian group. We will give briefly some basic definitions about Kleinian 
groups.\\

Let $M\ddot{o}b(\mathbb{S}^{n})$ denote the group of M\"obius
transformations of the n-sphere 
$\mathbb{S}^{n}=\mathbb{R}^{n}\cup\{\infty\}$, i.e. conformal diffeomorphisms of
$\mathbb{S}^{n}$ with respect to the standard metric. For a discrete group 
$G\subset M\ddot{o}b(\mathbb{S}^{n})$ {\it the discontinuity set} $\Omega(G)$ is
defined as follows
$$
\Omega (G)=\{x\in\mathbb{S}^{n}: \mbox{the point}\hspace{.2cm}
x\hspace{.2cm} \mbox{possesses a neighbourhood} 
\hspace {.2cm}U(x)\hspace{.2cm}\mbox{such that}
$$
$$
\hspace{1cm}
U(x)\cap g(U(x))\hspace{.2cm} 
\mbox{is empty for all but finite elements}\hspace{.2cm} g\in G\}.
$$
\vskip .3cm
The complement 
$\mathbb{S}^{n}-\Omega(G)=\Lambda(G)$ is called the {\it  limit set} (see 
\cite{kap1}).\\ 

Both $\Omega(G)$ and $\Lambda(G)$ are $G$-invariant, $\Omega(G)$ is
open, hence $\Lambda(G)$ is compact.\\

A subgroup $G\subset M\ddot{o}b(\mathbb{S}^{n})$ is called {\it  Kleinian} if 
$\Omega (G)$ is not empty. We will be concerned with very specific
Kleinian groups of Schottky type.\\

We recall that a conformal map $\psi$ on $\mathbb{S}^{n}$ can be extended in a
natural way to  the hyperbolic space $\mathbb{H}^{n+1}$,
such that $\psi|_{\mathbb{H}^{n+1}}$ is an orientation-preserving isometry with
respect to the Poincar\'{e} metric. Hence we can identify 
the group $M\ddot{o}b(\mathbb{S}^{n})$ with the group of orientation
preserving isometries of hyperbolic $(n+1)$-space $\mathbb{H}^{n+1}$. This
allows us to define the limit set of a Kleinian group through sequences.\\ 
 
A point $x$ is a {\it  limit point} for the Kleinian group
$G$, if there exist a point $z\in\mathbb{S}^{n}$ and a sequence $\{g_{m}\}$ 
of {\it distinct elements} of $G$, with $g_{m}(z)\rightarrow x$. 
The set of limit points is $\Lambda (G)$ (see \cite{maskit} section II.D).\\

One way to illustrate the action of a Kleinian group $G$ is to draw a picture 
of $\Omega(G)/G$. For this purpose a fundamental domain is very helpful. Roughly 
speaking, it contains one point from each equivalence class in $\Omega(G)$ (see 
\cite{kap2} pages  78-79, \cite{maskit} pages  29-30).\\

\begin{definition}
A fundamental domain $D$ for a Kleinian group $G$ is a 
co\-di\-men\-sion- zero piecewise-smooth submanifold (subpolyhedron) of
$\Omega(G)$ satisfying the following
\begin{enumerate}
\item $\bigcup_{g\in G} g(Cl_{\Omega(G)} D)=\Omega$ ($Cl$
  denotes closure).
\item $g(int(D))\cap int(D)=\emptyset$ for all
$g\in G-\{e\}$ ($int$ denotes the interior).
\item The boundary of $D$ in $\Omega (G)$ is a piecewise-smooth
(polyhedron) submanifold in $\Omega (G)$, divided into a union
of smooth submanifolds (convex polygons) which are called {\rm
faces}. For
each face $S$, there is a corresponding face $F$ and an element $g=g_{SF}\in
G-\{e\}$ such that $gS=F$ ($g$ is called a  face-pairing
transformation); $g_{SF}=g_{FS}^{-1}$. 
\item Only finitely many translates of $D$ meet any compact subset of
$\Omega (G)$.
 \end{enumerate} 
\end{definition}

\begin{theo}(\cite{kap2}, \cite{maskit}) Let 
$D^{*}=\overline{D}\cap\Omega/\sim_{G}$ denote the orbit space with
the quotient topology. Then $D^{*}$ is homeomorphic to $\Omega/G$.
\end{theo}

\section{The Construction}

The main idea of this construction is to use the symmetry of the spinning
process to find a ``packing'' (i.e. a cover) of an embedded $\mathbb{S}^{2}$ 
in $\mathbb{S}^{4}$ consisting of closed round balls
of dimension 4, such that the group $\Gamma$ generated by inversions in 
their boundaries (spheres of dimension 3) is Kleinian, and its limit set is a wild 
sphere of dimension two.\\

\begin{definition}
Let $X\in\mathbb{S}^{n}$. We will say that $E=\cup^{m}_{i=1} B_{i}$ is a {\it packing} 
for $X$ if this is contained
in the interior of $E$, where $B_{i}$ is a closed ball of dimension $n$ for $i=1,\ldots, m$. 
\end{definition}

\begin{definition}
Let $A$ be a spinnable knotted arc in $\mathbb{R}^{3}_{+}$. A {\it semi-pearl solid  necklace} 
subordinate to it, is a collection of consecutive closed round 4-balls $B^{1},\ldots,B^{n}$ such that 
\begin{enumerate}
\item The end-points of $A$ are the centers of  $B^{1}$ and $B^{n}$ respectively.
\item The arc $A$ is totally contained in $\cup_{i=1}^{n}B^{i}$.
\item Two consecutive balls are orthogonal; otherwise $B^{i}\cap
  B^{j}=\emptyset$,  $j\neq i+1$.
\item The segment of $A$ lying in the interior of each ball is unknotted.
\end{enumerate}
Each ball $B^{i}$ is called a {\it solid pearl}. Its boundary is a 3-sphere $\Sigma^{i}$ called 
a {\it pearl}. A {\it semi-pearl necklace} subordinate to $A$ is $\cup_{i=1}^{n}\Sigma^{i}$.
 
\end{definition}

Next, we will define a pearl-necklace $Spin(T)$ subordinate to the 2-knot $Spin(A)$, with
all the requirements needed for the group, generated by reflections on each pearl, to be
Kleinian.\\

\begin{definition}
Let $A\subset\mathbb{R}^{3}_{+}$ be a spinnable knotted arc. A pearl-necklace $Spin(T)$
subordinate to the tame knot $Spin(A)\subset\mathbb{S}^{4}$ is constructed in the
following way
\begin{enumerate}
\item Let $\overline{T}$ be a semi-pearl necklace subordinate to $A$ consisting of the pearls 
$\Sigma^{0},\Sigma^{1},\ldots,\Sigma^{l+1}$. Consider the subset  
$T=\{\Sigma^{1},\ldots,\Sigma^{l}\}$.
 Now, in $Spin(A)$ we will select six isometric
  copies $R_{2\pi i/6}(A)$ of $A$, called 
$A_{i}$, $i=1,\ldots,6$, in such a way that each $A_{i}$ has 
 subordinate a isometric copy $R_{2\pi i/6}(T)$ of $T$, 
denoted by $T_{i}$, $i=1,\ldots,6$. We will require that 
$\Sigma^{k}_{i}\in T_{i}$ is orthogonal to the corresponding $\Sigma^{k}_{i+1}\in T_{i+1}$.
\item At each pole of the knot $Spin(A)$ we set a pearl $\Sigma_{m}$, $m=1,2$, 
orthogonal to the pearls of the next and previous levels, $\Sigma^{1}_{i}$ and 
$\Sigma^{l}_{i}$, $i=1,\ldots,6$, respectively; 
such that $\Sigma_{1}\cap(\cap_{s=i}^{i+1}\Sigma^{1}_{s})\neq\emptyset$ and 
$\Sigma_{2}\cap(\cap_{s=i}^{i+1}\Sigma^{l}_{s})\neq\emptyset$
\item At each intersection point (see proposition 3.4) of two consecutive pearls
$\Sigma^{k}_{i}$, $\Sigma^{k+1}_{i}$ in $T_{i}$ and the corresponding
$\Sigma^{k}_{i+1}$, $\Sigma^{k+1}_{i+1}$ in $T_{i+1}$, we set a
pearl $P^{k}_{i}$ which is orthogonal to these four pearls and does
not intersect any other, i.e. we  
require that $P^{k}_{i}\cap\Sigma^{r}_{s}=\emptyset$ for 
$r\neq k,k+1$, $s\neq i,i+1$ and $P^{k}_{i}\cap\Sigma_{m}=\emptyset$, for $m=1,2$.
\item The intersection $B^{k}_{i}\cap Spin(A)$ is an unknotted disk, where $B^{k}_{i}$
is the solid pearl whose boundary is $\Sigma^{k}_{i}$.
\end{enumerate}
\end{definition}

Let $Spin(T)=\{\Sigma_{1},\ldots,\Sigma_{n}\}$ be a pearl-necklace. We define the {\it
  filling} of $Spin(T)$ as $|Spin(T)|=\cup_{i=1}^{n}B_{i}$,
where $B_{i}$ is the round closed 4-ball whose boundary $\partial B_{i}$
is the pearl $\Sigma_{i}$.\\

Geometrically, the above definition means that when we rotate $A$ and $T$ with respect to
$\mathbb{R}^{2}$, we obtain an infinite number of pearls covering
$Spin(A)$. We shall select a finite number of them keeping a ``symmetry'', i.e.
we will choose $T$ such that when we spin a pearl $\Sigma^{k}\in T$, we can 
select six of $\{R_{\theta}(\Sigma^{k}):0\leq\theta\leq 2\pi\}$ in
such a way that their centers form a regular hexagon and 
adjacent pearls are orthogonal (two pearls are orthogonal if the square of the
distance between their centers is equal to the sum of the squares of 
their radii). In other words, we will choose six $\mathbb{R}^{3}_{+}$ (six pages
of the open book decomposition of $\mathbb{R}^{4}$). As a consequence, in $Spin(A)$ we will have six 
preferential meridians (one for each page). Each meridian $A_{i}$ ($1\leq i\leq 6$), is
a copy of $A$ and has a semi pearl-necklace $T_{i}$ ($1\leq i\leq 6$) subordinate to it which is
an isometric copy of $T$. The pearls belong to $T_{i}$ will be denoted by $\Sigma^{k}_{i}$, where the
superscript $k=1,\ldots,l$, indicates its ``latitude'' and the subscript $i=1,\ldots,6$,
indicates its ``meridian'' (see Figure 2 and compare \cite{maskit} page 208).\\

\begin{figure}[tbh]
\centerline{\epsfxsize=2in \epsfbox{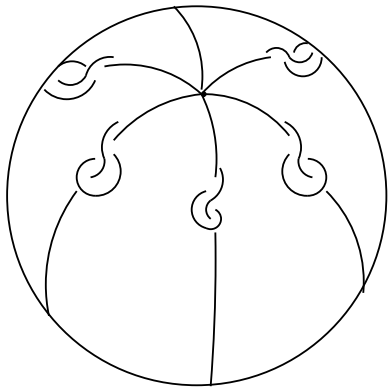}}
\caption{\sl Spinning the arc $A$ with six preferential meridians.}
\label{F2}
\end{figure}

At each pole of the knot $Spin(A)$ we will set a pearl $\Sigma_{m}$ ($m=1,2$)
orthogonal to the other six of the previous or next level respectively,
such that there is no hole among them (see Figure 3), i.e. 
$\Sigma_{m}\cap\Sigma^{j}_{i}\cap\Sigma^{j}_{i+1}\neq\emptyset$ ($i\in
\mathbb{Z}/6\mathbb{Z}$), for the 3-tuples of indexes  $(m=1,j=1,i)$
and $(m=2,j=l,i)$. By standard arguments of Euclidean geometry this sphere always exists.\\
  
\begin{figure}[tbh]
\centerline{\epsfxsize=1.5in \epsfbox{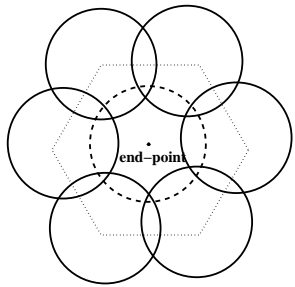}}
\caption{\sl A pearl set in an end-point of the arc $A$.}
\label{F3}
\end{figure}

At this point, we have chosen a finite number of pearls. However, we have not proved that $Spin(A)$
is totally covered by them.\\

\begin{prop}
The pearls $\Sigma^{k}_{i}$, 
$k=1,\ldots,l$, $i=1,\ldots,6$ and the two pearls $\Sigma_{m}$,
($m=1,2$) at the poles, totally cover a knot isotopic to $Spin(A)$.  
\end{prop}

{\bf {\it Proof.}}\\

Firstly, we will verify that the intersection of the pearls $\Sigma_{i}^{k}$,
$\Sigma_{i}^{k+1}\in T_{i}$ and the corresponding 
$\Sigma_{i+1}^{k}$, $\Sigma_{i+1}^{k+1}\in T_{i+1}$ is not empty.\\

Let $r$ be the radius of the pearls $\Sigma_{i}^{k}$ and $\Sigma_{i+1}^{k}$
with centers $c^{k}_{i}$ and $c^{k}_{i+1}$ respectively. Let $R$ be the radius
of the pearls $\Sigma_{i}^{k+1}$ and $\Sigma_{i+1}^{k+1}$ with centers 
 $c^{k+1}_{i}$ and $c^{k+1}_{i+1}$ respectively (see Figure 4). We know that  
two pearls are orthogonal if the square of the
distance between their centers is equal to the sum of the squares of 
their radii. Hence
$$
d^{2}(c^{k}_{i},c^{k}_{i+1})=2r^{2}
$$
$$
d^{2}(c^{k+1}_{i},c^{k+1}_{i+1})=2R^{2}
$$
$$
d^{2}(c^{k}_{i},c^{k+1}_{i})=d^{2}(c^{k}_{i+1},c^{k+1}_{i+1})=r^{2}+R^{2}
$$
this implies that
$$
d(c^{k}_{i},c^{k+1}_{i+1})=d(c^{k}_{i+1},c^{k+1}_{i})=r+R
$$
so that the intersection of these four pearls is a point.\\

\begin{figure}[tbh]
\centerline{\epsfxsize=1.2in \epsfbox{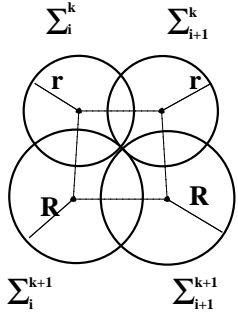}}
\caption{\sl Two consecutive pearls of $T_{i}$ and the correspondings of 
$T_{i+1}$.}
\label{F24}
\end{figure} 

Now consider the filling of each $T_{i}$. We shall prove
that there exists a knot isotopic to $Spin(A)$ which is totally contained
in $\mathcal{E}:=(\cup_{i=1}^{6}|T_{i}|)\cup_{m=1}^{2}B_{m}$, where $B_{m}$
is the closed 4-ball whose boundary is the pearl $\Sigma_{m}$.\\

Let $\omega:I\rightarrow \mathbb{R}^{3}_{+}$ be a parametrization of the arc $A$ by $t\in [0,1]$.
Let $\Pi_{t}\subset\mathbb{R}^{4}$ be an affine plane parallel to the $zw$-plane. 
Notice that $Spin(\omega(t))\subset \Pi_{t}$ is a circle for $t\in (0,1)$ and
a point for $t=0,1$.\\

Let $t=\epsilon_{1}>0$ be the smallest $t$ for which $\omega(\epsilon_{1})\in
B^{1}_{1}\cap\Sigma_{1}$ (remember that $B^{k}_{i}$ is the 4-ball such that
$\partial B^{k}_{i}=\Sigma^{k}_{i}$). Let $t=\epsilon_{2}>0$ be the smallest
$t$ that satisfies $\omega(\epsilon_{2})\in B^{l}_{1}\cap\Sigma_{2}$. For 
$t\in [\epsilon_{1},\epsilon_{2}]$, we have that $\omega(t)\in B^{k}_{1}$ for some
index $k$. Then $S_{t}:=\Pi_{t}\cap_{i=1}^{6}|\Sigma^{k}_{i}|$ is a union of
six disks with the property that adjacent disks are either overlapped or tangent.
Observe that $Spin(\omega(t))$ may be not contained in $S_{t}$ (see Figures 5 and 6).\\

\begin{figure}[tbh]
\centerline{\epsfxsize=1.5in \epsfbox{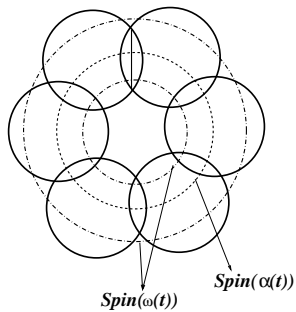}}
\caption{\sl $S_{t}$ is formed by six overlapping disks.}
\label{F25}
\end{figure}

\begin{figure}[tbh]
\centerline{\epsfxsize=1.5in \epsfbox{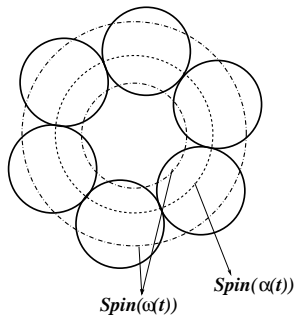}}
\caption{\sl $S_{t}$ is formed by six tangent disks.}
\label{F26}
\end{figure} 

By an isotopy of $\Pi_{t}$, we can send $Spin(\omega(t))$ to a circle $Spin(\alpha(t))$, which passes
through either the middle point of each chord joining the two intersection points of each
overlap or the points of tangency of adjacent circles (see Figures 5 and 6). Indeed, this
isotopy $\phi_{t}$ can be constructed radially from a function $\psi_{t}$ whose graph
appears in Figure 7 (both cases). Thus $\phi_{t}(s,x)=s(\psi_{t}(x))+(1-s)x$ is a
stable isotopy, i.e. is the identity in the complement of a closed set.\\

\begin{figure}[tbh]
\centerline{\epsfxsize=3.5in \epsfbox{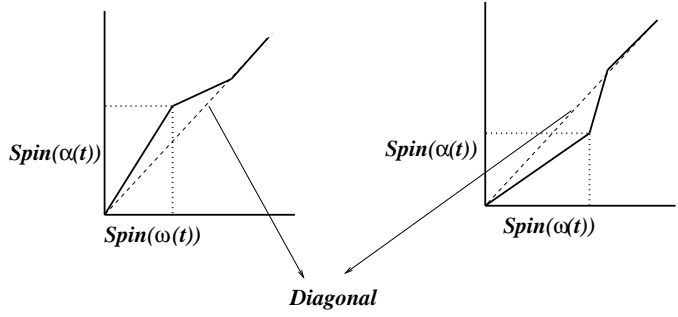}}
\caption{\sl Radial isotopy at level $\omega(t)$.}
\label{F27}
\end{figure}

For the pearls at the poles, we have that $Spin(\omega(t))\subset B_{1}$ for
$t\in [0,\epsilon_{1}]$ and $Spin(\omega(t))\subset B_{2}$ for
$t\in[\epsilon_{2},1]$. We can transform, by an isotopy, 
$\omega(t)$ $t\in [0,\epsilon_{1}]\cup [\epsilon_{2},1]$ in two arcs contained in the
interior of the respective balls with the condition that their end-points coincide
with $\alpha(\epsilon_{1})$, $\omega(0)$ and $\alpha(\epsilon_{2})$, $\omega(1)$, respectively.\\

By the above, we can define a function such that in each level $\omega(t)$ is the
previous isotopy. This function depends of the parameter of the isotopy on each level 
and $t$. Since it is continuous with respect to each variable, it is continuous.
Notice that on each level $\omega(t)$, we have that the corresponding isotopy is the
identity in the complement of some disk. Hence we can conclude that this function is the
identity in the complement of a closed ball.\\

We can extend this function to an isotopy defined on $\mathbb{S}^{4}$ (see \cite{palais})
that sends $Spin(\omega(t))$ to $Spin(\alpha(t))$.\\

Therefore, an isotopic knot to $Spin(A)$ is totally covered by the pearls 
$\Sigma^{k}_{i}$, $k=1,\ldots,l$, $i=1,\ldots,6$ and the
two pearls at the poles. $\blacksquare$\\

The intersection  of four pearls $\Sigma_{i}^{k}$,
$\Sigma_{i}^{k+1}\in T_{i}$ and   
$\Sigma_{i+1}^{k}$, $\Sigma_{i+1}^{k+1}\in T_{i+1}$ is a single point that will be denoted by
$p^{k}_{i}$ (see Figure 8). We centered at $p^{k}_{i}$ 
a pearl $P^{k}_{i}$ orthogonal to these four pearls such that it does not overlap to any
other. Notice that this sphere always exists and its
construction uses standard Euclidean geometry.\\

\begin{figure}[tbh]
\centerline{\epsfxsize=1.5in \epsfbox{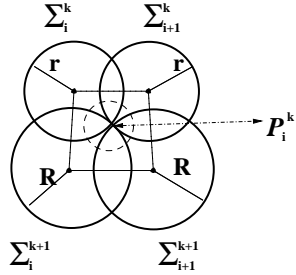}}
\caption{\sl The dotted sphere is the pearl $P^{k}_{i}$.}
\label{F28}
\end{figure}
 
Hence, $Spin(T)$ consists of the pearls $\Sigma^{k}_{i}$, $k=1,\ldots,l$, $i=1,\ldots,6$,
the two pearls $\Sigma_{m}$, $m=1,2$, at the poles and the pearls $P^{k}_{i}$. We will say
that $Spin(A)$ is the {\it template} of $Spin(T)$.\\

Consider the group $\Gamma$ generated by reflections through each pearl. To guarantee that  
the group $\Gamma$ is Kleinian we will use the Poincar\'e Polyhedron Theorem.
This theorem establishes conditions for the group to be discrete. In practice
these conditions are very hard to be verify, but in our case all of them are 
satisfied automatically from the construction (see \cite{kap2}, 
\cite{maskit}, \cite{epstein}).\\

This theorem also gives us a presentation for the group $\Gamma$. Suppose that the 
pearl-necklace $Spin(T)$ is formed by the pearls $\Sigma_{j}$, $(j=1,\ldots,n)$ and we
denote by $I_{j}$ the reflection with respect to $\Sigma_{j}$. Since the dihedral angles
between the faces $F_{i}$, $F_{j}$ are $\frac{\pi}{n_{ij}}$, where $n_{ij}$ is
either $2$ if the faces are adjacent or $0$ in other case. Therefore ,
we have the following presentation of $\Gamma$
$$
\Gamma=<I_{j}, j=1,\dots,n|\hspace{.2cm}(I_{j})^{2}=1,\hspace{.2cm}(I_{i}I_{j})^{n_{ij}}=1>
$$

\begin{prop}
The group $\Gamma$ generated by reflections through each pearl, is Kleinian.
\end{prop}

{\bf {\it Proof.}} By the Poincar\'e Polyhedron Theorem, we have that
$\Gamma$ is discrete and its fundamental domain is 
$\mathbb{S}^{4}-|Spin(T)|$. Therefore it is Kleinian. $\blacksquare$\\

The first question to appear is if there exists a pearl-necklace
$Spin(T)$ for some knot $Spin(A)$.
In the next theorem we will exhibit a semi-necklace $T$ subordinate to
an embedded of the trefoil
arc $A$, satisfying all the requirements of the definition 3.3.\\

\begin{theo}
There exists an embedding of the trefoil arc $A$ in $\mathbb{R}^{3}_{+}$ that admits a 
semi-necklace satisfying all the requirements of the definition 3.3. 
\end{theo}

{\bf {\it Proof.}} By proposition 3.4 it follows that if we have
constructed a pearl-necklace $Spin(T)$ subordinate
to the knot $Spin(A)$, it is
always possible to find a knot isotopic to  $Spin(A)$ such that it is totally contained in the interior of
$Spin(T)$. The group $\Gamma$ is defined through the pearl-necklace, this means that the 
pearl-necklace is more fundamental for our purpose than the knot
itself. This allows us to consider the trefoil arc $A$ as a
polygonal arc (see Figure 9) obtained joining the centers  $c_{k}$ of the pearls
$\Sigma^{k}_{1}$ whose coordinates appear in the next table.\\

\begin{figure}[tbh]
\centerline{\epsfxsize=3in \epsfbox{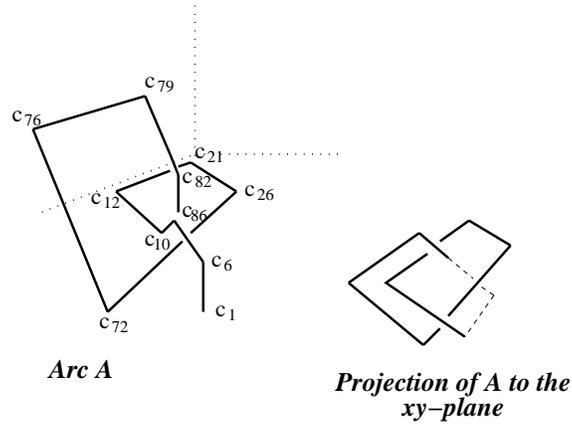}}
\caption{\sl A drawing of the trefoil-arc $A$.}
\label{F29}
\end{figure}

Observe that the pearl $\Sigma^{k}_{i}$ and the corresponding rotated $\Sigma^{k}_{i+1}$
are orthogonal if and only if their radii are equal to the $z$-coordinate divided by 
$\sqrt{2}$.\\

\begin{tabular}{|l|l|l|}\hline
$\Sigma_{1} $&$ C_{1}$=(2426.06421, 2296.89168, .75966995, 0) &
$R_{1}$=.537167778\\\hline

$\Sigma_{2} $&$ C_{2}$=(2426.06421, 2296.89168, 2.835126878,0) &$
R_{2}$=2.004737441\\\hline

$\Sigma_{3} $&$ C_{3}$=(2426.06421, 2296.89168, 10.58083755,0) &$
R_{3}$=7.481781981\\\hline

$\Sigma_{4} $&$ C_{4}$=(2426.06421, 2296.89168, 39.48822332,0) &$
R_{4}$=27.92239049\\\hline

$\Sigma_{5} $&$ C_{5}$=(2426.06421, 2296.89168, 147.3720558,0) &$
R_{5}$=104.20778\\\hline

$\Sigma_{6} $&$ C_{6}$=(2426.06421, 2296.89168, 550,0) &$
R_{6}$=388.9087297\\\hline

$\Sigma_{7} $&$ C_{7}$=(2426.06421, 1746.89168, 550,0)&$
R_{7}$=388.9087297\\\hline

$\Sigma_{8} $&$
C_{8}$=(2426.06421, 1196.89168, 550,0)&$R_{8}$=388.9087297\\\hline

$\Sigma_{9} $&$ C_{9}$=(2426.06421, 740, 400,0)&$
R_{9}$=282.8427125\\\hline

$\Sigma_{10} $&$ C_{10}$=(2226.56071, 597.879, 200,0)&$
R_{10}$=141.4213562\\\hline

$\Sigma_{11} $&$ C_{11}$=(2126.56071, 397.879, 258.5786444,0)&$
R_{11}$=182.8427129\\\hline

$\Sigma_{12} $&$ C_{12}$=(2026.56071, 197.879, 200,0)&$
R_{12}$=141.4213562\\\hline

$\Sigma_{13} $&$ C_{13}$=(1826.56071, 197.879, 200,0)&$
R_{13}$=141.4213562\\\hline

$\Sigma_{14} $&$ C_{14}$=(1626.56071, 197.879, 200,0)&$
R_{14}$=141.4213562\\\hline

$\Sigma_{15} $&$ C_{15}$=(1426.56071, 197.879, 200,0)&$
R_{15}$=141.4213562\\\hline

$\Sigma_{16} $&$ C_{16}$=(1226.56071, 197.879, 200,0)&$
R_{16}$=141.4213562\\\hline

$\Sigma_{17} $&$ C_{17}$=(1026.56071, 197.879, 200,0)&$
R_{17}$=141.4213562\\\hline

$\Sigma_{18} $&$ C_{18}$=(826.56071, 197.879, 200,0)&$
R_{18}$=141.4213562\\\hline

$\Sigma_{19} $&$ C_{19}$=(626.56071, 197.879, 200,0)&$
R_{19}$=141.4213562\\\hline

$\Sigma_{20} $&$ C_{20}$=(426.56071, 197.879, 200,0)&$
R_{20}$=141.4213562\\\hline

$\Sigma_{21} $&$ C_{21}$=(426.56071, 390.842826, 186.6225781,0)&$
R_{21}$=131.9620905\\\hline

$\Sigma_{22} $&$ C_{22}$=(426.56071, 556.1388641, 150,0)&$
R_{22}$=106.0660172\\\hline

$\Sigma_{23} $&$ C_{23}$=(426.56071, 695.063304, 130,0)&$
R_{23}$=91.92388155\\\hline

$\Sigma_{24} $&$ C_{24}$=(390, 804.611533, 105,0)&$
R_{24}$=74.24621202\\\hline
\end{tabular}

\begin{tabular}{|l|l|l|}\hline
$\Sigma_{25} $&$ C_{25}$=(420.7470362, 905.0088369, 105,0)&$
R_{25}$=74.24621202\\\hline

$\Sigma_{26} $&$ C_{26}$=(518.54878, 907.99193, 91.98807,0)&$
R_{26}$=65.0453\\\hline

$\Sigma_{27} $&$ C_{27}$=(610.5368498, 907.99193, 91.98807,0)&$
R_{27}$=65.0453\\\hline

$\Sigma_{28} $&$ C_{28}$=(702.5249198, 907.99193, 91.98807,0)&$
R_{28}$=65.0453\\\hline

$\Sigma_{29} $&$ C_{29}$=(794.5129898, 907.99193, 91.98807,0)&$
R_{29}$=65.0453\\\hline

$\Sigma_{30} $&$ C_{30}$=(886.5010598, 907.99193, 91.98807,0)&$
R_{30}$=65.0453\\\hline

$\Sigma_{31} $&$ C_{31}$=(978.4891298, 907.99193, 91.98807,0)&$
R_{31}$=65.0453\\\hline

$\Sigma_{32} $&$ C_{32}$=(1070.4772, 907.99193, 91.98807,0)&$
R_{32}$=65.0453\\\hline

$\Sigma_{33} $&$ C_{33}$=(1162.46527, 907.99193, 91.98807,0)&$
R_{33}$=65.0453\\\hline

$\Sigma_{34} $&$ C_{34}$=(1254.45334, 907.99193, 91.98807,0)&$
R_{34}$=65.0453\\\hline

$\Sigma_{35} $&$ C_{35}$=(1346.44141, 907.99193, 91.98807,0)&$
R_{35}$=65.0453\\\hline

$\Sigma_{36} $&$ C_{36}$=(1438.42948, 907.99193, 91.98807,0)&$
R_{36}$=65.0453\\\hline

$\Sigma_{37} $&$ C_{37}$=(1530.41755, 907.99193, 91.98807,0)&$
R_{37}$=65.0453\\\hline

$\Sigma_{38} $&$ C_{38}$=(1622.40562, 907.99193, 91.98807,0)&$
R_{38}$=65.0453\\\hline

$\Sigma_{39} $&$ C_{39}$=(1714.39369, 907.99193, 91.98807,0)&$
R_{39}$=65.0453\\\hline

$\Sigma_{40} $&$ C_{40}$=(1806.38176, 907.99193, 91.98807,0)&$
R_{40}$=65.0453\\\hline

$\Sigma_{41} $&$ C_{41}$=(1898.36983, 907.99193, 91.98807,0)&$
R_{41}$=65.0453\\\hline

$\Sigma_{42} $&$ C_{42}$=(1990.3579, 907.99193, 91.98807,0)&$
R_{42}$=65.0453\\\hline

$\Sigma_{43} $&$ C_{43}$=(2082.34597, 907.99193, 91.98807,0)&$
R_{43}$=65.0453\\\hline

$\Sigma_{44} $&$ C_{44}$=(2174.33404, 907.99193, 91.98807,0)&$
R_{44}$=65.0453\\\hline

$\Sigma_{45} $&$ C_{45}$=(2266.32211, 907.99193, 91.98807,0)&$
R_{45}$=65.0453\\\hline

$\Sigma_{46} $&$ C_{46}$=(2358.31018, 907.99193, 91.98807,0)&$
R_{46}$=65.0453\\\hline

$\Sigma_{47} $&$ C_{47}$=(2450.29825, 907.99193, 91.98807,0)&$
R_{47}$=65.0453\\\hline

$\Sigma_{48} $&$ C_{48}$=(2542.28632, 907.99193, 91.98807,0)&$
R_{48}$=65.0453\\\hline

$\Sigma_{49} $&$ C_{49}$=(2634.27439, 907.99193, 91.98807,0)&$
R_{49}$=65.0453\\\hline

$\Sigma_{50} $&$ C_{50}$=(2726.26246, 907.99193, 91.98807,0)&$
R_{50}$=65.0453\\\hline

$\Sigma_{51} $&$ C_{51}$=(2818.25053, 907.99193, 91.98807,0)&$
R_{51}$=65.0453\\\hline

$\Sigma_{52} $&$ C_{52}$=(2910.2386, 907.99193, 91.98807,0)&$
R_{52}$=65.0453\\\hline

$\Sigma_{53} $&$ C_{53}$=(3002.22667, 907.99193, 91.98807,0)&$
R_{53}$=65.0453\\\hline

$\Sigma_{54} $&$ C_{54}$=(3094.21474, 907.99193, 91.98807,0)&$
R_{54}$=65.0453\\\hline

$\Sigma_{55} $&$ C_{55}$=(3186.20281, 907.99193, 91.98807,0)&$
R_{55}$=65.0453\\\hline

$\Sigma_{56} $&$ C_{56}$=(3278.19088, 907.99193, 91.98807,0)&$
R_{56}$=65.0453\\\hline

$\Sigma_{57} $&$ C_{57}$=(3370.17895, 907.99193, 91.98807,0)&$
R_{57}$=65.0453\\\hline

$\Sigma_{58} $&$ C_{58}$=(3462.16702, 907.99193, 91.98807,0)&$
R_{58}$=65.0453\\\hline

$\Sigma_{59} $&$ C_{59}$=(3554.15509, 907.99193, 91.98807,0)&$
R_{59}$=65.0453\\\hline

$\Sigma_{60} $&$ C_{60}$=(3646.14316, 907.99193, 91.98807,0)&$
R_{60}$=65.0453\\\hline

$\Sigma_{61} $&$ C_{61}$=(3738.13123, 907.99193, 91.98807,0)&$
R_{61}$=65.0453\\\hline
\end{tabular}

\begin{tabular}{|l|l|l|}\hline
$\Sigma_{62} $&$ C_{62}$=(3830.1193, 907.99193, 91.98807,0)&$
R_{62}$=65.0453\\\hline

$\Sigma_{63} $&$ C_{63}$=(3922.10737, 907.99193, 91.98807,0)&$
R_{63}$=65.0453\\\hline

$\Sigma_{64} $&$ C_{64}$=(4014.09544, 907.99193, 91.98807,0)&$
R_{64}$=65.0453\\\hline

$\Sigma_{65} $&$ C_{65}$=(4106.08351, 907.99193, 91.98807,0)&$
R_{65}$=65.0453\\\hline

$\Sigma_{66} $&$ C_{66}$=(4198.07158, 907.99193, 91.98807,0)&$
R_{66}$=65.0453\\\hline

$\Sigma_{67} $&$ C_{67}$=(4290.05965, 907.99193, 91.98807,0)&$
R_{67}$=65.0453\\\hline

$\Sigma_{68} $&$ C_{68}$=(4382.04772, 907.99193, 91.98807,0)&$
R_{68}$=65.0453\\\hline

$\Sigma_{69} $&$ C_{69}$=(4474.03579, 907.99193, 91.98807,0)&$
R_{69}$=65.0453\\\hline

$\Sigma_{70} $&$ C_{70}$=(4566.02386, 907.99193, 91.98807,0)&$ 
R_{70}$=65.0453\\\hline

$\Sigma_{71} $&$ C_{71}$=(4650.89687, 943.4652579, 91.98807,0)&$ 
R_{71}$=65.0453\\\hline

$\Sigma_{72} $&$ C_{72}$=(4750, 914.7460538, 120,0)&$ R_{72}$=84.85281374\\\hline
$\Sigma_{73} $&$ C_{73}$=(4750, 770.0255859, 203.776124,0)&$
 R_{73}$=144.0914791\\\hline

$\Sigma_{74} $&$ C_{74}$=(4750, 522.386491, 451.4162296,0)&$
R_{74}$=319.19947\\\hline

$\Sigma_{75} $&$ C_{75}$=(4570, 0, 930.302113,0)&$R_{75}$=657.8229327\\\hline
$\Sigma_{76} $&$ C_{76}$=(3611.943, -100, 1000,0)&$ R_{76}$=707.1067812\\\hline
$\Sigma_{77} $&$ C_{77}$=(2611.943, -100, 1000,0)&$ R_{77}$=707.1067812\\\hline
$\Sigma_{78} $&$ C_{78}$=(1611.943, -100, 1000,0)&$ R_{78}$=707.1067812\\\hline
$\Sigma_{79} $&$ C_{79}$=(1096.504432, 756.930197, 1000,0)&$ R_{79}$=707.1067812\\\hline
$\Sigma_{80} $&$ C_{80}$=(922.2268, 600, 300,0)&$ R_{80}$=212.1320344\\\hline
$\Sigma_{81} $&$ C_{81}$=(750, 500, 163.335525,0)&$ R_{81}$=115.4956573\\\hline
$\Sigma_{82} $&$ C_{82}$=(750, 500, 43.765626,0)&$ R_{82}$=30.94697\\\hline
$\Sigma_{83} $&$ C_{83}$=(750, 500, 11.72696421,0)&$ R_{83}$=8.29221591\\\hline
$\Sigma_{84} $&$ C_{84}$=(750, 500, 3.1422119,0)&$ R_{84}$=2.221879342\\\hline
$\Sigma_{85} $&$ C_{85}$=(750, 500, .841953143,0)&$ R_{85}$=.595350776\\\hline      
\end{tabular}

\begin{rem}
The coordinates of the centers of $\Sigma_{76}$, $\Sigma_{77}$ and $\Sigma_{78}$
are rational numbers and their radii are equal to $\frac{1000}{\sqrt{2}}$. We
obtained the rest of centers and radii using the equations
$$
d^{2}(C_{k-1},C_{k})=R_{k-1}^{2}+R_{k}^{2}
$$
and
$$
d^{2}(C_{k+1},C_{k})=R_{k+1}^{2}+R_{k}^{2}
$$
Hence, we conclude that all centers and radii of the pearls
belong to a finite algebraic extension of the rational numbers.
\end{rem}




Let $\Gamma$ be the group generated by reflections $I_{j}$, through the pearl 
$\Sigma_{j}$ ($j=1,\ldots,n$) of the necklace $Spin(T)$ formed by $n$ pearls. 
Then $\Gamma$ is a conformal Kleinian group.\\

\section{Geometric Description of the Limit Set}

Let $A$ be a spinnable knotted arc in $\mathbb{S}^{3}$. Consider the 2-knot 
$Spin(A)\subset\mathbb{S}^{4}$ and take a pearl-necklace $Spin(T)$ subordinate to $Spin(A)$ consisting
on $n$ pearls.\\

Let $\Gamma$ be the group generated by reflections $I_{j}$ through the pearl $\Sigma_{j}\in Spin(T)$.
The natural question is: What is its limit set? Recall that 
to find the limit set of $\Gamma$, we need to find all the accumulation points of orbits. To
do that we are going to consider all the possible sequences of elements of   
$\Gamma$. We will do this in steps:

\begin{enumerate}
\item First step: Reflecting with respect to each $\Sigma_{j}$ 
($j=1,2,\ldots,n$), a copy of the exterior of $Spin(T)$ is mapped within it. 
At the end we obtain a new knot
$Spin(A_{1})$, which is in turn isotopic to the connected sum of $n+1$ copies 
of $Spin(A)$ and it is totally covered by  $n(n-2)$ pearls (packing) called
$E(T_{1})$.\\

Notice that there exists an isotopy of $\mathbb{S}^{4}$ such that the knot 
$Spin(A\# A)$ is sent to $Spin(A)\# Spin(A)$. Actually, this remains true 
for the connected sum of any couple of knotted arcs. 
Therefore $Spin(A_{1})$ is isotopic to spin of the connected sum of $n+1$
copies of $A$.\\

We have that
$|Spin(T)|=|E(T_{1})|$. In fact, each pearl of $Spin(T)$ lies in 
$E(T_{1})$.  Hence $|E(T_{1})|$ is a closed neighbourhood of 
$Spin(A)$. To clarify the above, see Figure 10 for a simpler case, 
i.e. for an unknotted necklace.\\
 
\begin{figure}[tbh]
\centerline{\epsfxsize=2in \epsfbox{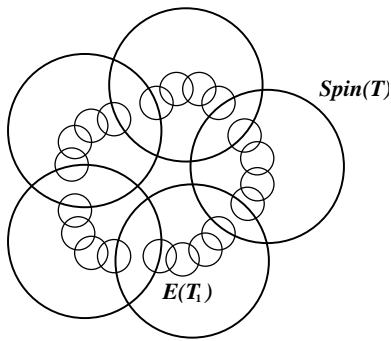}}
\caption{\sl An unknotted necklace and the first iteration.}
\label{F32}
\end{figure}

\begin{cla}
$|Spin(T_{1})|=V$ is isotopic to a closed tubular neighbourhood of 
$Spin(A)$.
\end{cla}

{\it Proof.} Let $N$ be a closed tubular neighbourhood of $Spin(A)$
 with the condition that $N\subset Int (V)$. 
Since $A$ is a spinnable arc, it follows that $Spin(A)$ is smooth. Given $p\in Spin(A)$,
consider the tangent plane $T_{p}Spin(A)$ of 
$Spin(A)$ at $p$. Let $\Pi^{2}(p)\subset \mathbb{S}^{4}$ be a 2-sphere totally geodesic with
respect to the spherical metric (i.e. radius 1) that passes through $p$ and intersects transversally
$T_{p}Spin(A)$. Thus $\Pi^{2}(p)$ cuts each solid pearl $B_{i}\in |Spin(T)|$ that contains
$p$ in a disk $D_{i}$. Then $D_{p}=\cup D_{i}$ is a star-shaped set
with respect to $p$. This neighbourhood is contained in a closed disk $\cal{B}$$_{R_{p}}(p)$,
where the radius is $R_{p}=\sup\{d(x,p):x\in D_{p}\} + \epsilon$. Notice that
$N_{p}=\Pi^{2}(p)\cap N\subset D_{p}$ (see Figure 11).\\ 

\begin{figure}[tbh]
\centerline{\epsfxsize=1.8in \epsfbox{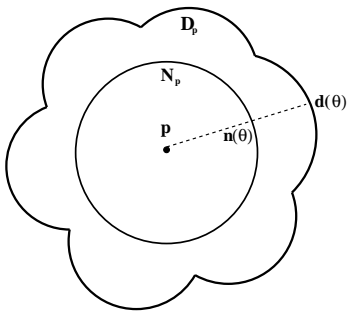}}
\caption{\sl A start-shaped neighbourhood of $p$.}
\label{F33}
\end{figure}

Hence for each point $p\in Spin(A)$, we have found a neighbourhood $D_{p}\subset\Pi^{2}(p)$ 
of it which is star-shaped with respect to $p$ and retracts onto $N_{p}$. This retraction can
be constructed in the following way. One draws a ray $r_{\theta}$ going from $p$ with angle $\theta$.
Let $n(\theta)$ be the intersection point of $r_{\theta}$ with $N_{p}$ and let $d(\theta)$ be
the intersection point of $r_{\theta}$ with $D_{p}$ (see Figure 13). We can send the segment
$d(\theta)$ to the segment $n(\theta)$ through the radial isotopy 
$\phi_{p}(t,x)=t\psi_{p}(x)+(1-t)x$, where $\psi_{p}(s)$ is the unique polygonal function whose
graph appears in Figure 12. Observe that this function is the identity beyond a distance
$R_{p}$ from $p$.\\

\begin{figure}[tbh]
\centerline{\epsfxsize=1.8in \epsfbox{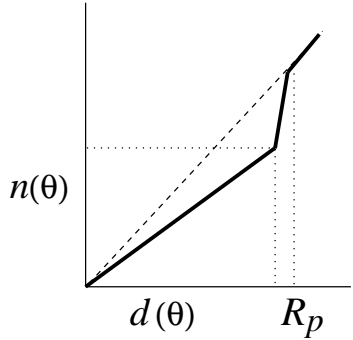}}
\caption{\sl The ray $d(\theta)$ is sent to the ray $n(\theta)$
by a radial isotopy.}
\label{F34}
\end{figure}
 
Therefore we have an isotopy defined on $\Pi^{2}(p)$ that transforms $D_{p}$ to $N_{p}$ and is the
identity outside of some closed disk $\cal{B}$$_{R_{p}}(p)$. Since $\Pi^{2}(p)$ depends
continuously  of $p$, we have an isotopy that sends $V$ to $N$. By \cite{palais} we
can extend this isotopy to $\mathbb{S}^{4}$. This isotopy is the identity outside of a closed 
tubular neighborhood $Spin(A)\times\mathbb{S}^{1}$, where the radius of $\mathbb{S}^{1}$ is
$\sup\{R_{p}\}$. $\blacksquare$\\  

\item Second step: If we consider the action of elements of $\Gamma$
on $E(T_{1})$, we obtain a new knot $Spin(A_{2})$ totally covered by a packing consisting of
$n(n^{2}-2n+7)$ pearls, called $E(T_{2})$. 
The knot $Spin(A_{2})$ is isotopic to the connected
sum of $n^{2}+1$ copies of $Spin(A)$. By the above observation, it
follows that it is also 
isotopic to the Spin of the connected sum of $n^{2}+1$ copies of
$A$.\\

Let $V_{1}=E(T_{2})-Spin(T)$. Then $|V_{1}|$ is connected and
is a closed neighbourhood of the 2-knot $Spin(P_{1})$, which is in turn isotopic to the
connected sum of $2n+1$ copies of $Spin(A)$. By the above claim,
$|V_{1}|$ is isotopic to a closed tubular 
neighbourhood of $Spin(P_{1})$. Notice that $|V_{1}|\subset |V|$ (see Figure 13).\\

\begin{figure}[tbh]
\centerline{\epsfxsize=2in \epsfbox{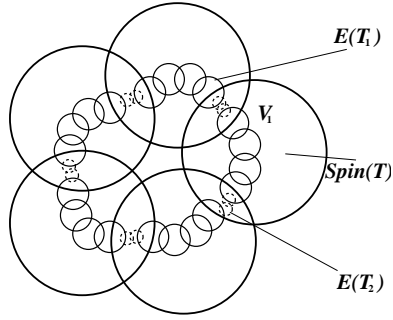}}
\caption{\sl The dotted pearls form $|V_{1}|$.}
\label{F35}
\end{figure}

\item $k^{th}$-Step: The action of elements of $\Gamma$ on $E(T_{k-1})$ 
determines a tame knot $Spin(A_{k})$, which is in turn isotopic to the connected 
sum of $n[\frac{(n-1)^{k}-1}{n-2}]+1$
copies of $Spin(A)$ and is also isotopic to the Spin of the connected sum of 
$n[\frac{(n-1)^{k}-1}{n-2}]+1$ copies of $A$.\\

Let $V_{k-1}=E(T_{k})-E(T_{k-2})$. Thus $|V_{k-1}|$ is connected and is
a closed neighbourhood of the knot $Spin(P_{k-1})$ which is in turn isotopic to
the connected sum of $n[\frac{(n-1)^{k-1}-1}{n-2}]+n(n-3)^{k-2}+1$. This
neighbourhood consists of $2n(n-3)^{k}$ pearls and is isotopic to a
closed tubular neighbourhood of $Spin(P_{k-1})$. By construction, 
$|V_{k-1}|\subset |V_{k-2}|$.

\end{enumerate}

Let $x\in\cap_{k=1}^{\infty} |V_{k}|$. We shall prove that $x$ is a
 limit point. Indeed, there exists a sequence of closed balls
 $\{B_{m}\}$ with $B_{m}\subset |V_{m}|$ such that $x\in B_{m}$ for
 each $m$. We can find a $z\in\mathbb{S}^{4}-Spin(T)$
 and a sequence $\{w_{m}\}$ of distinct elements of $\Gamma$, such
 that $w_{m}(z)\in B_{m}$. Since $diam(B_{m})\rightarrow 0$ it follows
 that $w_{m}(z)$ converges to $x$. The other inclusion clearly holds. 
Therefore, the limit set is given by
$$
\Lambda (\Gamma , A)=\varprojlim_{k} |V_{k}|=\bigcap_{k=1}^{\infty} |V_{k}|.
$$

\begin{theo} The limit set $\Lambda(\Gamma , A)$ is isotopic to $Spin(\Lambda)$,
where $\Lambda$ is a wild arc in the sense of \cite{kap1}, \cite{maskit}, and is
contained in each page ($\mathbb{R}^{3}_{+}$) of the open book decomposition of 
$\mathbb{R}^{4}$. 
\end{theo}

{\bf {\it Proof.}}
Let $A$ be a spinnable knotted arc. Construct the 2-knot $Spin(A)$ and take the
necklace of $n$-pearls $Spin(T)$ subordinate to $Spin(A)$.\\

Now consider a semi-pearl necklace $C$ consisting of $n$ consecutive orthogonal round
2-spheres that cover completely to $A$, in which its end-points are the 
centers of the first pearl, $\Sigma_{1}$, and the last one, $\Sigma_{n}$.  
Construct $Spin(C)=\cup_{0\leq\theta\leq 2\pi}R_{\theta}(C)$ (see section 2).\\

\begin{cla}
$|Spin(T)|$ is isotopic to $|Spin(C)|$.
\end{cla}

Indeed, we have already proved that $|Spin(T)|$ is isotopic to a closed tubular
neighbourhood of the knot $Spin(A)$. By the same argument, $|Spin(C)|$ 
is isotopic to a closed tubular neighbourhood of $Spin(A)$. Now two closed tubular
neighbourhoods of $Spin(A)$ are isotopic (\cite{hirsch}). This proves the claim. $\blacksquare$\\

In the first step of the reflecting process applied to $Spin(T)$, 
we get a packing $E(T_{1})$, of 
$Spin(A_{1})$ formed by pearls. Now, for the case of the semi-necklace $C$, we join the end-points 
of $A$ by an unknotted curve $L$ obtaining a knot $K$ (see Figure 14).\\

\begin{figure}[tbh]
\centerline{\epsfxsize=1.5in \epsfbox{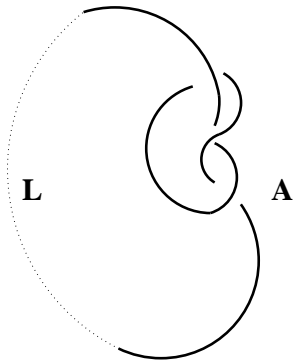}}
\caption{\sl The dotted curve $L$ joining the end-points of the arc $A$.}
\label{F36}
\end{figure}

We complete the semi-pearl necklace $C$ for the knot $K$, with pearls $Z_{s}$, $s=1,\ldots,r$, keeping
the same conditions on consecutive pearls. This new necklace is called $Z$ (see Figure 15).\\
 
\begin{figure}[tbh]
\centerline{\epsfxsize=1.5in \epsfbox{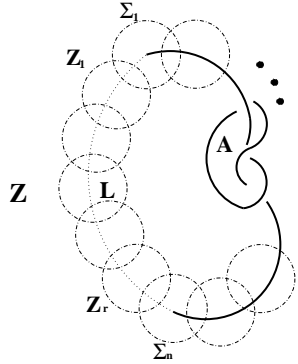}}
\caption{\sl The pearl-necklace $Z$ subordinate to the knot $K$.}
\label{F37}
\end{figure}

Now, we reflect only with respect to each pearl $\Sigma_{i}$, $i=1,\ldots,n$, of $C$. Then we obtain 
a new knot $K^{1}$ isotopic to the connected sum of $n+1$ copies of $K$. 
To return $K^{1}$ to an arc, we remove the unknotted curve joining the image of the end-points 
of $A$ under the corresponding reflections. This new arc is called
$A^{1}$ and is totally covered by a set of pearls $C_{1}$. Observe
that $C\subset C_{1}$ (see Figure 16).\\
 
\begin{figure}[tbh]
\centerline{\epsfxsize=1.5in \epsfbox{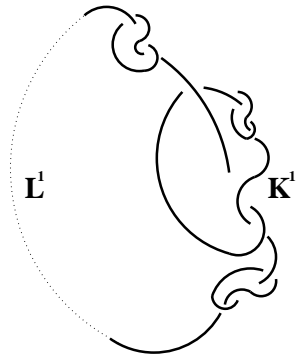}}
\caption{\sl The knot $K^{1}$ and the line $L^{1}$.}
\label{F38}
\end{figure}

Notice that
$Spin(A^{1})$ is isotopic to $Spin(A_{1})$ and    
$Spin(C_{1})$ is a packing for it. Thus, $Spin(|C_{1}|)=Spin(|C|)$
(where $|C|$ is defined as $|Spin(T)|$) is a clo\-sed neigh\-bour\-hood of 
$Spin(A)$. \\

In the second step for the necklace $Spin(T)$, we get a packing $E(T_{2})$ of the knot $Spin(A_{2})$. 
For the semi-necklace $C$, we join again the end-points of the arc $A^{1}$
by an unknotted curve, $L_{1}$, forming again the knot $K^{1}$ in such a way that when we complete
the semi-pearl necklace $C_{1}$, we add  the pearls $Z_{s}$ $s=1,\ldots,r$, obtaining  
the necklace $Z$. We can assume that the 
end-points of the arc $A^{1}$ coincide with the centers of the pearls $\Sigma_{1}$ and 
$\Sigma_{n}\in C\subset C_{1}$, respectively (see Figure 17).\\
 
\begin{figure}[tbh]
\centerline{\epsfxsize=1.5in \epsfbox{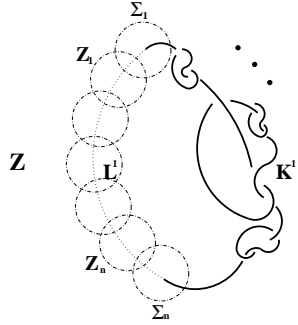}}
\caption{\sl The pearl-necklace $Z$ subordinate to the knot $K^{1}$.}
\label{F39}
\end{figure}

Now reflecting only with respect to each pearl
of the semi-necklace $C$, i.e. with respect to the pearls $\Sigma_{i}$
$i=1,\ldots,n$, we get as in the previous step, 
the packing $C_{2}$ of the new arc
$A^{2}$, which is in turn isotopic to the connected sum of $2n+1$ copies of $K$ minus an 
unknotted curve $L_{2}$. When we spin $C_{2}$ and $A^{2}$, we obtain the packing $Spin(C_{2})$ 
of $Spin(A^{2})$. Define $W_{1}=Spin(C_{2})-Spin(C)\cup\{I_{j}(Z):\hspace{.2cm}j=1,\ldots,n\}$. 
Then $|W_{1}|$ is a closed neighbourhood of a 2-knot $Spin(Q^{2})$, which is isotopic to the connected sum
of $2n+1$ copies of $Spin(A)$. So $|V_{1}|$ and $|W_{1}|$
are closed  neighbourhoods of isotopic 2-knots. Using the same arguments 
of claim 4.1 and the standard fact that any locally flat embedding
of $\mathbb{S}^{2}$ in $\mathbb{S}^{4}$ has trivial normal bundle, it follows that
two closed tubular neighbourhoods of 
isotopic knots are isotopic, hence $|V_{1}|$ is isotopic to 
$|W_{1}|$ and the following diagram commutes
$$
\begin{CD}
(\mathbb{S}^{4},|V_{1}|)@>>>(\mathbb{S}^{4},|V|)\\
@V\sim VV@V\sim VV\\
(\mathbb{S}^{4},|W_{1}|)@>>>(\mathbb{S}^{4},|W|),
\end{CD}
$$      
where the row maps are inclusions. Notice that this isotopy is stable,
i.e. is the identity on some open in $\mathbb{S}^{4}$,
and is orientation-preserving (see \cite{kirby}).\\

Inductively, for the 
$k^{th}$-step we obtain the packings $E(T_{k})$ of $Spin(A_{k})$
and $Spin(C_{k})$ of $Spin(A^{k})$. Where $Spin(A_{k})$ is obtained through the
reflecting process previously described. The arc $A^{k}$ is formed applying the
reflecting process to $C\subset Z$ subordinate to the knot
$K^{k}$ and removing an unknotted curve $L_{k}$. Then
$$
|V_{k}|=|E(T_{k+1})-E(T_{k-1})|
$$ 
and
$$
|W_{k}|=|Spin(C_{k+1})-Spin(C_{k-1})\cup\{I_{i_{l}i_{l-1}\cdots i_{1}}(Z):\hspace{.2cm}1\leq l\leq k\}|
$$
are closed neighborhoods of the knots $Spin(P_{k})$ and $Spin(Q^{k})$ respectively, which are in turn
isotopic to the connected sum of $n[\frac{(n-1)^{k}-1}{n-2}]+n(n-3)^{k-1}+1$ copies of $Spin(A)$. 
Hence, $|V_{k}|$ is isotopic to $|W_{k}|$ and the following diagram commutes
$$
\begin{CD}
(\mathbb{S}^{4},|V_{k}|)@>>>(\mathbb{S}^{4},|V_{k-1}|)\\
@V\sim VV@V\sim VV\\
(\mathbb{S}^{4},|W_{k}|)@>>>(\mathbb{S}^{4},|W_{k-1}|).
\end{CD}
$$          
Observe that this isotopy is stable and orientation-preserving.
Summarizing, we have the commutative diagram
$$
\begin{CD}
(\mathbb{S}^{4},|V|)@<<<(\mathbb{S}^{4},|V_{1}|)@<<<\cdots @<<<(\mathbb{S}^{4},|V_{k}|)\cdots\\
@V\sim VV@V\sim VV@V\sim VV@V\sim VV\\
(\mathbb{S}^{4},|W|)@<<<(\mathbb{S}^{4},|W_{1}|)@<<<\cdots
@<<<(\mathbb{S}^{4},|W_{k}|)\cdots\\
\end{CD}
$$
where the row maps are inclusions and the vertical arrows are
orientation-preserving stable isotopies.\\

The inverse limit in the first row of the above diagram is
$(\mathbb{S}^{4},\,\Lambda(\Gamma , A))$ and the inverse limit in the second row is
$(\mathbb{S}^{4},\, Spin(\varprojlim_{k}|W_{k}|))$. But  
$\varprojlim_{k}|W_{k}|$ is a wild arc de\-no\-ted by
$\Lambda(\Gamma)$ (see \cite{maskit}, \cite{kap2}), i.e. the
in\-ver\-se limit in the second row is 
$(\mathbb{S}^{4},\,Spin(\Lambda(\Gamma)))$.\\ 

By the universal property of the inverse limit, 
there exists a homeomorphism of $\mathbb{S}^{4}$ to $\mathbb{S}^{4}$ which sends 
$\Lambda(\Gamma , A))$ to $Spin(\Lambda(\Gamma))$. This homeomorphism is stable because it coincides with 
a stable homeomorphism on some open set (see \cite{kirby}) and is orientation-preserving. This implies that
it is isotopic to the identity (see \cite{kirby}).\\ 

Therefore, the knots
$Spin(\Lambda(\Gamma))$ and $\Lambda(\Gamma , A))$ are isotopic. This
proves the Theorem.$\blacksquare$

 \begin{coro}
The limit set $\Lambda(\Gamma , A)$ is homeomorphic to $\mathbb{S}^{2}$
\end{coro}

{\bf {\it Proof.}}
By the above theorem, we have that
$$
\Lambda(\Gamma , A)\cong Spin(\Lambda(\Gamma))\cong \mathbb{S}^{2}.
$$ 
\rightline{$\blacksquare$}

\begin{theo} Let $Spin(T)$ be a pearl-necklace subordinate to the non-trivial 
tame knot $Spin(A)$. Then  $\Lambda(\Gamma , A)$ is wildly embedded in
$\mathbb{S}^{4}$.
\end{theo}

{\bf {\it Proof.}}
The fundamental group of $\mathbb{S}^{4}-\Lambda(\Gamma , A)$ is isomorphic to the fundamental
group of the knot obtained joining the end-points of the arc $\Lambda$ by an unknotted curve 
(see \cite{rolfsen}). It is well-known that this fundamental group  has no finite
representation (see \cite{kap1}, \cite{maskit}).  $\blacksquare$ 

\begin{ex}
Let $Spin(T)$ be a pearl-necklace subordinate to
$Spin(A)$ where $A$ is the trefoil arc, $T_{2,3}$. Then
$$
\Pi_{1}(Spin(A))\cong\Pi_{1}(T_{2,3})=\{x,y\mid xyx=yxy\}
$$
hence
$$
\hspace{-5.5cm}\Pi_{1}(\mathbb{S}^{4}-\Lambda(\Gamma,A))=
\{x_{1},y_{1},\ldots,y_{n},\ldots \mid
$$
$$
\hspace{1cm}x_{1}y_{1}x_{1}=y_{1}x_{1}y_{1},
x_{1}y_{2}x_{1}=y_{2}x_{1}y_{2},\ldots,x_{1}y_{n}x_{1}=y_{n}x_{1}y_{n}
,\ldots \}
$$
$$
\cong(\cdots(\Pi_{1}(Spin(A))*_{\{x_{1}\}}\Pi_{1}(Spin(A)))*_{\{x_{1}\}}
\cdots *_{\{x_{1}\}}\Pi_{1}(Spin(A))*_{\{x_{1}\}} \cdots
$$
is infinitely generated with a infinite number of relations.
\end{ex}

\section{Hyperbolic Manifolds}

The action of $\Gamma$ can be extended to the hyperbolic space $\mathbb{H}^{5}$ and in this
case $\Gamma$ is a subgroup of Isom $\mathbb{H}^{5}$, which acts properly and discontinuously
on $\mathbb{D}^{5}=\mathbb{H}^{5}\cup\partial\mathbb{H}^{5}$. Its fundamental polyhedron is
${\cal{P}}=(\mathbb{H}^{5}\cup\partial\mathbb{H}^{5})-\widetilde{|Spin(T)|}$, where
$\widetilde{Spin(T)}$ is the natural extension of the pearl-necklace to $\mathbb{H}^{5}$. It
is a convex subset and has a finite number of sides, hence $\Gamma$ is geometrically finite
(see \cite{finites}).\\ 

The group $\Gamma$ acts properly and discontinuously on $\overline{\cal{P}}$, then the quotient
${\cal{M}}^{5}_{\Gamma}=(\mathbb{D}^{5}-\Lambda(\Gamma,A))/\Gamma\cong\overline{\cal{P}}$
(see Theorem 2.2) is a compact orbifold such that its interior is a non-compact hyperbolic manifold
of infinite volume and its compactification as a subset of $\mathbb{D}^{5}$ has boundary which
possesses a conformally flat structure given by the action.\\

For the Kleinian group $\Gamma$ acting on the pearl-necklace $Spin(T)$, its fundamental domain is
$D=\mathbb{S}^{4}-|Spin(T)|$. The group $\Gamma$ acts properly and discontinuously on $\overline{D}$,
hence $\overline{D}\cong\Omega(\Gamma)/\Gamma=(\mathbb{S}^{4}-\Lambda(\Gamma))/\Gamma$ is
an orientable, compact, conformally flat 4-orbifold with boundary. Its fundamental group
coincides with the fundamental group of the template of $Spin(T)$.\\

In the next section, we will describe $(\mathbb{S}^{4}-\Lambda(\Gamma))/\Gamma$ under 
the restriction that $Spin(A)$ is a fibered knot.\\

Consider now the index-two subgroup  $\widetilde{\Gamma}\subset\Gamma$ consisting of even words, i.e.
$\widetilde{\Gamma}$ is the orientation preserving index two subgroup of $\Gamma$. Its 
fundamental polyhedron is 
$\widetilde{\cal{P}}=(\mathbb{H}^{5}\cup\partial\mathbb{H}^{5}-\widetilde{|Spin(T)|})\cup
(\widetilde{B_{j}}-\widetilde{I_{j}(|Spin(T)-\Sigma_{j}|)})$, where tilde means the natural
extensions to the hyperbolic space of both the pearl-necklace and the corresponding reflection map.
Since $\widetilde{\cal{P}}\subset\mathbb{D}^{5}$ is a convex subset and has a finite number
of sides, it follows that $\widetilde{\Gamma}$ is geometrically finite.\\

Since $\widetilde{\Gamma}$ acts freely on its domain of discontinuity,
then the quotient space
${\cal{M}}^{5}_{\widetilde{\Gamma}}=(\mathbb{D}^{5}-\Lambda(\widetilde{\Gamma},A))/
\widetilde{\Gamma}\cong\widetilde{P}/_{\sim\widetilde{\Gamma}}$
is a compact, orientable manifold, such that $Int({\cal{M}}^{5}_{\widetilde{\Gamma}})$
is a non-compact, orientable hyperbolic manifold of infinite volume. This space as a subset of
$(\mathbb{D}^{5}-\Lambda(\widetilde{\Gamma}))$, 
has a boundary which 
possesses a natural conformally flat structure given by the action.\\

For the Kleinian group $\widetilde{\Gamma}$ acting on
$\mathbb{S}^{4}$, its fundamental domain is 
$\widetilde{D}=(\mathbb{S}^{4}-|Spin(T)|)\cup (B_{j}-I_{j}(|Spin(T)-\Sigma_{j}|))$.
Since $\widetilde{\Gamma}$ acts freely on $\Omega(\widetilde{\Gamma})$, we have
that $\Omega/\widetilde{\Gamma}\cong \overline{D}\cap\Omega/\sim_{\widetilde{\Gamma}}$ 
is a compact, orientable, conformally flat 4-manifold with boundary. Its fundamental group is the
fundamental group of the knot $Spin(A\# A)$.

\section{Fibration of $\mathbb{S}^{4}-\Lambda(\Gamma)$ over
 $\mathbb{S}^{1}$}

We recall that a mapping $f:E\rightarrow B$ is said to be a locally trivial fibration with 
fiber $F$ if each point of $B$ has a neighbourhood $U$ and a ``trivializing'' 
homeomorphism $h:f^{-1}(U)\rightarrow U\times F$ for which the following diagram 
commutes
$$
\xymatrix{
f^{-1}(U)\ar[d]_f\ar[r]^h& U\times F\ar[dl]^-{\txt{projection}}\\
U&                                                             
}
$$
$E$ and $B$ are known as the {\it total} and {\it base} spaces, respectively. 
Each set $f^{-1}(b)$ is called a {\it fiber} and is homeomorphic to $F$. We will 
be concerned with fibrations with base space $\mathbb{S}^{1}$.\\

\begin{definition}
A knot or link $L$ in $\mathbb{S}^{3}$ is {\it fibered} if there exists a 
locally trivial fibration $f:(\mathbb{S}^{3}-L)\rightarrow \mathbb{S}^{1}$. We 
require that $f$ be well-behaved near $L$. That is, each component $L_{i}$ is 
to have a neighbourhood framed as $\mathbb{D}^{2}\times\mathbb{S}^{1}$, with 
$L_{i}\cong \{0\}\times\mathbb{S}^{1}$, in such a way that the
restriction of $f$ 
to $(\mathbb{D}^{2}-\{0\})\times\mathbb{S}^{1}$ is the map into $\mathbb{S}^{1}$ 
given by $(x,y)\rightarrow \frac{y}{|y|}$.\\
\end{definition}
It follows that each $f^{-1}(x)\cup L$, $x\in\mathbb{S}^{1}$, is a 2-manifold
 with boundary $L$: in fact a Seifert surface for $L$ 
(see \cite{rolfsen}, page 323).\\

\begin{ex}
Let $\mathbb{S}^{3}_{\epsilon}\subset\mathbb{C}^{2}$ be the 3-sphere
centered at the origin of radius $\epsilon$. Let 
$V=\{(z_{1},z_{2})\in\mathbb{C}^{2}:z_{1}^{2}+z_{2}^{3}=0\}$. Then 
$\mathbb{S}^{3}_{\epsilon}\cap V=K$ is the right-handed trefoil knot and the map 
$F:\mathbb{S}^{3}_{\epsilon}-K\rightarrow\mathbb{S}^{1}$ given by 
$F(z_{1},z_{2})=\frac{z_{1}^{2}+z_{2}^{3}}{|z_{1}^{2}+z_{2}^{3}|}$ is a 
locally trivial fibration with fiber the punctured torus (see \cite{milnor} 
section 1, \cite{rolfsen} pages  327-333). 
\end{ex}
\vskip .3cm
\begin{lem}
(\cite{zeeman1})Let $A$ be a spinnable knotted arc. Suppose that the knot $K$, obtained from $A$ 
joining its end-points by an unknotted curve, fibers over the
circle with fiber the surface $S$. Then $Spin(A)$ fibers over the
circle with fiber $S_{\theta}$, an $\mathbb{S}^{1}$-family of surfaces
$\tilde{S}$ all glued onto a single meridian of
$\partial\mathbb{D}^{3}$ with longitude $\theta$. The
interior of $\tilde{S}$ is $S$ and its boundary is a  meridian of 
$\partial\mathbb{D}^{3}$ (see Figure 18). 
\end{lem}

{\bf {\it Proof.}}
The fibering of the complement of $K$ induces a fibering of 
$\mathbb{D}^{3}-\mathbb{D}^{1}$ by surfaces $\tilde{S}_{\theta}$,
$\theta\in\mathbb{S}^{1}$. The interior of $\tilde{S}_{\theta}$ is $S$
and its boundary is $\partial\tilde{S}_{\theta}$=$M^{1}_{\theta}$ the
meridian of $\partial\mathbb{D}^{3}$ with longitude $\theta$ (see
Figure 18).\\

\begin{figure}[tbh]
\centerline{\epsfxsize=1.5in \epsfbox{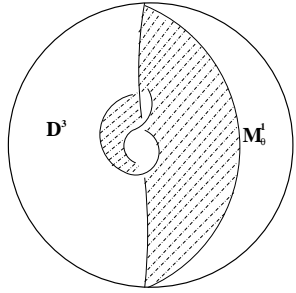}}
\caption{\sl Fibering of $\mathbb{D}^{3}-\mathbb{D}^{1}$.}
\label{F40}
\end{figure}

Recall that in the spinning process we multiply the interior of 
$\mathbb{D}^{3}$ by $\mathbb{S}^{1}$ and $\partial\mathbb{D}^{3}$
stays fixed. Hence, we get a fibering of $\mathbb{S}^{4}-Spin(A)$ by an 
$\mathbb{S}^{1}$-family of surfaces $\tilde{S}$ all glued onto the single
meridian $M^{1}_{\theta}$ of $\partial \mathbb{D}^{3}$ with longitude
$\theta$. \rightline{$\blacksquare$} 
\vskip .5cm
Henceforth, a fibered arc will mean that the knot obtained from it joining
its end-points by an unknotted curve, fibers over the circle.\\
 
\begin{lem}
Let $A$ be a spinnable fibered arc with fiber the surface $S$.
Let $Spin(T)$ be an $n$-pearl necklace subordinate to the tame knot $Spin(A)$.
Let $\Lambda(\Gamma , A)$ be the limit set. Then
$\Omega (\Gamma)/\Gamma$ fibers over the circle with fiber $S^{**}$, the
closure of the surface $S_{\theta}$ of the previous lemma.
\end{lem}

{\bf {\it Proof.}}
Let $\widetilde{P}:\mathbb{S}^{4}-Spin(A)\rightarrow\mathbb{S}^{1}$ be the
given fibration with fiber the 3-manifold $S_{\theta}$. Observe that 
$\widetilde{P}\mid_{\mathbb{S}^{4}-Spin(T)}\equiv P$ is a fibration with fiber
$S_{\theta}$.\\

As we know $\Omega (\Gamma)/\Gamma\cong D^{*}$. In our case 
$D^{*}=\overline{\mathbb{S}^{4}-Spin(T)}$, which fibers over the circle
with fiber the closure of $S_{\theta}$.$\blacksquare$\\

By the above, to describe $(\mathbb{S}^{4}-\Lambda(\Gamma, A))/\Gamma$ when the original knot 
is fibered, we just need to determine its monodromy. It coincides with the knot's monodromy.
Hence we have a complete description of $(\mathbb{S}^{4}-\Lambda(\Gamma, A))/\Gamma$.\\

\begin{lem}
Let $A$ be a spinnable fibered arc with fiber the surface $S$. 
Let $Spin(T)$ be an $n$-pearl necklace subordinate to the tame knot
$Spin(A)$. Let $\widetilde{\Gamma}$ be the orientation preserving
index two subgroup of $\Gamma$. 
Let $\Lambda(\widetilde{\Gamma} , A)$ be the limit set. Then
$\Omega (\widetilde{\Gamma})/\widetilde{\Gamma}$ fibers over the
circle with fiber $S^{*}$, which is homeomorphic to the connected sum along 
the boundary of the 3-manifold $S_{\theta}$ with itself.
\end{lem}

{\bf {\it Proof.}}
We can assume, up to isotopy, that the fiber $S$ cuts each pearl of
the semi-necklace corresponding to $A$, in arcs going from one
intersection point to another. Hence, we can assume that the fiber
$S_{\theta}$ cuts each pearl $\Sigma_{i}\in Spin(T)$ in disks $a_{i}$,
whose boundary is the intersection of $\Sigma_{i}$ with the adjacent
pearls.\\

When we reflect with respect to $\Sigma_{i}$ a copy of $S_{\theta}$,
called $S^{i}_{\theta}$, is mapped to the interior of $\Sigma_{i}$ and
it is joined to $S_{\theta}$ along the disk $a_{i}$. $\blacksquare$\\

Since $\widetilde{\Gamma}$ is a normal subgroup of $\Gamma$, it follows by Lemma 8.1.3 in 
\cite{thurston1} that $\widetilde{\Gamma}$ has the same limit set as $\Gamma$. Therefore
$\mathbb{S}^{4}-\Lambda(\Gamma,
A)=\mathbb{S}^{4}-\Lambda(\widetilde{\Gamma}, A)$.

\begin{theo}
Let $A$ be a non-trivial spinnable fibered arc. Let $Spin(T)$ be a pearl-necklace subordinate
to the fibered knot $Spin(A)$. Let $\Gamma$ be the group generated by
reflections through the pearls and let $\widetilde{\Gamma}$ be the
orientation preserving index two subgroup of $\Gamma$. Let 
$\Lambda(\Gamma , A)=\Lambda(\widetilde{\Gamma}, A)$ be the 
corresponding limit set. Then:
\begin{enumerate}
\item There exists a locally trivial fibration $\psi
:\mathbb{S}^{4}-\Lambda(\Gamma , A)\rightarrow\mathbb{S}^{1}$, where the
 fiber $\Sigma^{*}_{\theta}=\psi^{-1}(\theta)$ is an
$\mathbb{S}^{1}$-family of surfaces $\Sigma$ all glued onto a
meridian $\theta$, of $\partial\mathbb{D}^{3}$ (see Lemma 6.3). Where $\Sigma$
is an orientable infinite genus surface with one end.
\item  $\overline{\Sigma^{*}_{\theta}}-\Sigma^{*}_{\theta}=\Lambda(\Gamma , A)$.
\end{enumerate}
\end{theo}

{\bf {\it Proof.}}
We know that 
$\zeta:\Omega(\widetilde{\Gamma})\rightarrow\Omega(\widetilde{\Gamma})/\widetilde{\Gamma}$
is an infinite-fold covering. By the previous lemma, there exists a locally trivial 
fibration $\phi:\Omega(\widetilde{\Gamma})/\widetilde{\Gamma}\rightarrow\mathbb{S}^{1}$ with fiber 
$S^{*}$.\\

Then $\psi=\phi\circ\zeta:\Omega(\widetilde{\Gamma})\rightarrow\mathbb{S}^{1}$ is a 
locally trivial fibration. The fiber is $\Gamma(S^{*})$, i.e. the orbit of 
the fiber.\\

We now give another proof. As we know from Theorem 4.1, 
the knot $\Lambda(\Gamma , A)$ is
isotopic to the knot $Spin(\Lambda(\Gamma))$, where $\Lambda(\Gamma)$
is a wild arc. Since $A$ is fibered, so is $\Lambda(\Gamma)$. In this
case the fiber, $\Sigma$, is an orientable infinite genus surface with
one end. Hence $Spin(\Lambda(\Gamma))$ fibers over the circle with fiber
$\Sigma^{*}_{\theta}$, an $\mathbb{S}^{1}$-family of surfaces $\Sigma$ all 
glued onto a meridian $\theta$, of $\partial\mathbb{D}^{3}$ (see Lemma 6.3).\\

The first part of the theorem has been proved. For the second part,
observe that the closure of the fiber is the closure of the
$\mathbb{S}^{1}$-family of surfaces $\Sigma$, i.e. is the closure of
an $\mathbb{S}^{1}$-family of ends. As we can see in the Figure 20, each end
has as boundary the wild arc $\Lambda(\Gamma)$. Hence the closure of the fiber 
is exactly the limit set. Therefore 
$\overline{\Sigma^{*}_{\theta}}-\Sigma^{*}_{\theta}=\Lambda(\Gamma , A)$. $\blacksquare$

\begin{rem}
\begin{enumerate}
\item This theorem can be generalized to  fibered links.
\item This theorem gives an open book decomposition of
$\mathbb{S}^{4}-\Lambda(\Gamma , A)$, where the ``binding'' is the wild
knot $\Lambda(\Gamma , A)$, and each ``page'', $\Sigma^{*}$, is a
3-manifold which fibers over $\mathbb{S}^{1}$ and is the
$\mathbb{S}^{1}$-family of surfaces $\Sigma$ all glued onto a
meridian $\theta$, of $\mathbb{D}^{3}$ (see Lemma 6.3). Here $\Sigma$
is an orientable infinite genus surface with one end.\\
Indeed, this decomposition can be viewed in the following way. For the above 
theorem, $\mathbb{S}^{4}-\Lambda(\Gamma , A)$ is $\Sigma^{*}\times [0,1]$ modulo the 
identification of the top with the bottom through an identifying 
homeomorphism. Consider $\overline{\Sigma^{*}}\times [0,1]$ and identify
the top with the bottom. This is equivalent to keep  
$\partial\overline{\Sigma^{*}}$ fixs
and to spin $\Sigma^{*}\times \{0\}$ with respect to
$\partial\overline{\Sigma^{*}}$ until 
glue it  with $\Sigma^{*}\times \{1\}$. 
Removing $\partial\overline{\Sigma^{*}}$ we obtain the open book decomposition. 

\end{enumerate}
\end{rem}
\section{Monodromy}

Let $Spin(A)$ be a non-trivial fibered tame knot and let $S$ be the fiber. Since 
$\mathbb{S}^{4}-Spin(A)$ fibers over the circle, we know that
$\mathbb{S}^{4}-Spin(A)$ is a mapping torus equal to $S\times [0,1]$ modulo  
an identifying homeomorphism $\psi:S\rightarrow S$ that glues $S\times\{0\}$ to
$S\times\{1\}$. This homeomorphism induces a homomorphism
$$
\psi_{\#}:\Pi_{1}(S)\rightarrow\Pi_{1}(S)
$$
called {\it the monodromy of the fibration}.\\

Another way to understand the monodromy is through  the 
{\it first return Poincar\'e map}, defined as follows. Let $M$ be
connected, compact  
manifold and let $f_{t}$ be a flow that possesses a transversal 
section $\eta$. It follows that if $x\in\eta$ then there exists a
continuous function $t(x)>0$ such 
that $f_{t}\in\eta$. We may define the first return Poincar\'e map  
$F:\eta\rightarrow\eta$ as $F(x)=f_{t(x)}(x)$. This map is a diffeomorphism 
and induces a homomorphism of $\Pi_{1}$ 
called {\it the monodromy} (see \cite{verjovsky}, chapter 5).\\

For the manifold $\mathbb{S}^{4}-Spin(A)$, the flow that defines the first 
return Poincar\'e map $\Phi$ is the flow that cuts transversally each page of 
its open book decomposition.\\

Consider a pearl-necklace $Spin(T)$ subordinate to $Spin(A)$. As we have observed
during the reflecting process, $Spin(A)$ and $S$ are copied in each
reflection. So the flow $\Phi$ is also copied. Hence, the 
Poincar\'e map can be extended in each step, giving us in the end a
homeomorphism $\psi:\Sigma^{*}_{\theta}\rightarrow \Sigma^{*}_{\theta}$ that identifies 
$\Sigma^{*}_{\theta}\times\{0\}$ with $\Sigma^{*}_{\theta}\times\{1\}$, and induces
the monodromy of the wild knot.\\

From the above, if we know the monodromy of the knot $Spin(A)$ then we know
the monodromy of the wild knot $\Lambda (\Gamma ,A)$.\\

By the long exact sequence associated to a fibration, we have
\begin{equation}
0\rightarrow \Pi_{1}(\Sigma^{*}_{\theta})\rightarrow
\Pi_{1}(\mathbb{S}^{4}-\Lambda (\Gamma ,A))\stackrel{\Psi}{\overset{\longleftarrow}{\rightarrow}}
\mathbb{Z}\rightarrow 0, \tag{1}
\end{equation}
which has a homomorphism section 
$\Psi:\mathbb{Z}\rightarrow(\mathbb{S}^{4}-\Lambda (\Gamma ,A))$. 
Therefore (1) splits. As a consequence 
$\Pi_{1}(\mathbb{S}^{4}-\Lambda (\Gamma ,A))$ is
the semi-direct product of $\mathbb{Z}$ with $\Pi_{1}(\Sigma^{*}_{\theta})$.\\

\begin{ex}
Let $A$ be the trefoil arc. Consider the knot $Spin(A)$. Then the
fiber $\Sigma^{*}_{\theta}$ $(\theta\in\mathbb{S}^{1})$ is an
$\mathbb{S}^{1}$-family of punctured torus all glued onto a single meridian 
(see previous section). The
fundamental group of $\Sigma^{*}_{\theta}$ is the free group in two generators, $a$ and
$b$. Since $\Pi_{1}(Spin(A))\cong\Pi_{1}($Trefoil knot), it follows that the 
monodromy maps, in both cases, coincide. That is, $\psi_{\#}$ sends 
$a\mapsto b^{-1}$ and $b\mapsto ab$. Its order is six up to an outer 
automorphism (See \cite{rolfsen} pages 330-333).\\

The monodromy in the limit $\psi_{\#}:\Pi_{1}(\Sigma^{*}_{\theta})\rightarrow 
\Pi_{1}(\Sigma^{*}_{\theta})$ is given by $a_{i}\mapsto b^{-1}_{i}$ and 
$b_{i}\mapsto a_{i}b_{i}$, where
$\Pi_{1}(\Sigma^{*}_{\theta})=\{a_{i},b_{i}\}$. So
$$
\begin{aligned}
\Pi_{1}(\mathbb{S}^{4}-\Lambda (\Gamma ,A))&\cong
\Pi_{1}(\mathbb{S}^{1})\ltimes_{\psi_{\#}}\Pi_{1}(\Sigma^{*}_{\theta})\\
&=\{a_{i},b_{i},c:a_{i}*c=b_{i}^{-1},\hspace{.2cm}b_{i}*c=a_{i}b_{i}\}\\
&=\{a_{i},c:c^{-1}a_{i}^{-1}c=a_{i}c^{-1}a_{i}\}\\
&=\{a_{i},c:c=a_{i}ca_{i}c^{-1}a_{i}^{-1}\}\\
&=\{a_{i},c:c=ca_{i}ca_{i}c^{-1}a_{i}^{-1}c^{-1}\}\\
&=\{a_{i},c:c=ca_{i}c^{-1}c^{2}a_{i}c^{-1}c^{-1}ca_{i}^{-1}c^{-1}\}.\\
\end{aligned}
$$
Let $\alpha_{i}=ca_{i}c^{-1}$;\\

$$
\begin{aligned}
\hspace{.9cm}&=\{\alpha_{i},c:c=\alpha_{i}c\alpha_{i}c^{-1}\alpha_{i}^{-1}\}\\
&=\{\alpha_{i},c:c\alpha_{i}c=\alpha_{i}c\alpha_{i}\}
\end{aligned}
$$
This gives another method for computing the fundamental group of a wild 2-knot
whose complement fibers over the circle.
\end{ex}

\begin{coro}Let $Spin(T)$ be a pearl-necklace whose template is a
non-trivial tame fibered knot $Spin(A)$. Then
$\Pi_{1}(\Omega(\Gamma)/\Gamma)\cong\mathbb{Z}\ltimes_{\psi_{\#}}\Pi_{1}(\Sigma^{*}_{\theta})$.  
\end{coro}

\section{Kleinian Groups and Twistor Spaces}

In this section we will lift the action of the group $\Gamma$ on $\mathbb{S}^{4}$
to its twistorial space, which is complex projective 3-space
$P^{3}_{\mathbb{C}}$. We refer to  \cite{alberto1} and \cite{alberto2} for details.\\

Let us now recall briefly the {\it twistor fibration} of $\mathbb{S}^{4}$,  also known as
{\it the Calabi-Penrose fibration}
$\pi:P^{3}_{\mathbb{C}}\rightarrow\mathbb{S}^{4}$ (see \cite{alberto2}). There are
several equivalent ways to construct this fibration. A geometric way to describe it is
by thinking of $\mathbb{S}^{4}$ as being the quaternionic projective line 
${P}^{1}_{\cal{H}}$, of {\it right} quaternionic lines in the quaternionic plane
${\cal{H}}^{2}$ (regarded as a 2-dimensional right $\cal{H}$-module). That is, for 
$q:=(q_{1},q_{2})\in\cal{H}^{2}$ ($q\neq (0,0)$) the right quaternionic line passing
through $q$ is the linear space
$R_{q}:=\{(q_{1}\lambda ,q_{2}\lambda )|\lambda\in {\cal{H}}\}$.
We can identify ${\cal{H}}^{2}$ with $\mathbb{C}^{4}$ via the 
$\mathbb{R}$-linear map given by $(q_{1},q_{2})\mapsto (z_{1},z_{2},z_{3},z_{4})$, where
$q_{1}=z_{1}+z_{2}${\bf j}=$x_{1}+x_{2}${\bf i}+$x_{3}${\bf j}+$x_{4}${\bf k} and 
$q_{2}=z_{3}+z_{4}${\bf j}=$y_{1}+y_{2}${\bf i}+$y_{3}${\bf j}+$y_{4}${\bf k}. In this notation 
{\bf i}, {\bf j}, {\bf k} denote the standard quaternionic units, 
$z_{1}=x_{1}+x_{2}${\bf i}, $z_{2}=x_{3}+x_{4}${\bf i}, $z_{3}=y_{1}+y_{2}${\bf i} and
$z_{4}=y_{3}+y_{4}${\bf i}.\\

Under this identification each right quaternionic line is invariant under {\it right} multiplication
by {\bf i}. Hence such a line is canonically isomorphic to $\mathbb{C}^{2}$. If we think of
$P^{3}_{\mathbb{C}}$ as being the space of complex lines in $\mathbb{C}^{4}$, then there
is an obvious map $\pi:P^{3}_{\mathbb{C}}\rightarrow\mathbb{S}^{4}$, whose fiber over a point
$H\in {P}^{1}_{\cal{H}}$ is the space of complex lines in the given right quaternionic line
$H\cong\mathbb{C}^{2}$; thus the fiber is $P^{1}_{\mathbb{C}}$.\\

The group $Conf_{+}(\mathbb{S}^{4})$ of orientation preserving conformal automorphisms of 
$\mathbb{S}^{4}$ is isomorphic to $PSL(2,{\cal{H}})$, the projectivization of the group
$2\times 2$, invertible, quaternionic matrices. This is naturally a subgroup of 
$PSL(4,\mathbb{C})$, since every quaternion corresponds to a couple of complex numbers. 
Hence $Conf_{+}(\mathbb{S}^{4})$ has a canonical lifting to a group of holomorphic
transformations of $P^{3}_{\mathbb{C}}$, carrying twistor lines into twistor lines.

\begin{definition}
(\cite{alberto1}) By a {\rm twistor Kleinian group} we mean a discrete subgroup $G$ of  
$Aut_{hol}(P^{3}_{\mathbb{C}})$ of holomorphic automorhisms, which acts on 
$P^{3}_{\mathbb{C}}$ with non-empty region of discontinuity $\Omega(\Gamma)$ and 
which is a lifting of a conformal Kleinian group acting on $\mathbb{S}^{4}$.
\end{definition}

\begin{rem}
There is no  ``good'' general definition of the discontinuity set
$\Omega$ for general groups, hence an appropiate definition must be
given in each case (see \cite{kulkarni1}.
We are considering the definition 1.4 of \cite{alberto1}, in which $\Omega(G)$ is an open 
$G$-invariant set and $G$
acts properly and  discontinuously on $\Omega(G)$. The space $\Omega(G)/G$ has the quotient topology, 
and the map  $\pi:\Omega\rightarrow\Omega/G$ is continuous and open.
\end{rem}

It has been proved in \cite{alberto1} that if  $G\subset Conf_{+}(\mathbb{S}^{4})$ is a discrete subgroup
acting on $\mathbb{S}^{4}$ with limit set $\Lambda$, then its canonical lifting 
$\widetilde{Conf_{+}}(\mathbb{S}^{4})$ acts on $P^{3}_{\mathbb{C}}$ with limit set 
$\widetilde{\Lambda}=\pi^{-1}(\Lambda)$ , thus $\widetilde{\Lambda}$ is a fibered bundle over  
$\Lambda$ with fiber $\mathbb{S}^{2}$. In \cite{alberto1}, is also proved that if we restrict the 
twistor bundle to a proper subset of  $\mathbb{S}^{4}$.\\

We consider the Kleinian group  $\Gamma$ such that its limit set is $\mathbb{S}^{2}$
wildly embedded on $\mathbb{S}^{4}$. Then 

\begin{theo}
There exists a $\mathbb{S}^{2}\times\mathbb{S}^{2}$ wildly embedded in the twistor space 
$P^{3}_{\mathbb{C}}$ dynamically defined, i.e. it is the limit set of
a complex Kleinian group $\Gamma\subset Aut_{hol}(P^{3}_{\mathbb{C}})$. $\blacksquare$
\end{theo}


\begin{thebibliography}{99}
\bibitem{artin} E. Artin, \emph{Zur Isotopie zweidimensionalen
  Fl\"achen im $R_{4}$} Abh. Math. Sem. Univ. Hamburg
(1926), 174-177.
\bibitem{finites} B.H. Bowditch. \emph{Geometrical Finiteness for
Hyperbolic Groups}. Journal of Functional Analysis 113 (1993),
  245-317.
\bibitem{dugundji} J. Dugundji. \emph{Topology}. Allyn and Bacon, Inc.
 1966.
\bibitem{epstein} D. B. A. Epstein, C. Petronio. \emph{An exposition of 
Poincare's polyhedron theorem}.
Enseignement Mathematique 40, 1994, 113-170.
\bibitem{fox1} R. H. Fox. \emph{A Quick Trip Through Knot
Theory}. Topology of 3-Manifolds and Related Topics. Prentice-Hall,
Inc., 1962.
\bibitem{GLT} M. Gromov, H. B. Lawson, W. Thurston. \emph{Hyperbolic
4-manifolds and conformally flat 3-manifolds}. Publ. Math. I.H.E.S. Vol. 68
(1988), 27-45.
\bibitem{goldman} W. Goldman. \emph{Conformally Flat Manifolds
with Nilpotent Holonomy and the Uniformization Problem for
3-Manifolds}. Transactions of the American Mathematical Society
Vol. 278 No. 2, 573-583.
\bibitem{hirsch} W. Hirsch. \emph{Smooth Regular Neighbourhoods}. Annals
of Mathematics Vol. 76, No.3 (1962), 524-530.
\bibitem{kap1} M. Kapovich. \emph{Topological Aspects of Kleinian
Groups in Several Dimensions}. Preprint (1988).
\bibitem{kap2} M. Kapovich. \emph{Hyperbolic Manifolds and
Discrete Groups}. Progress in Mathematics, Birkhauser, 2001.
\bibitem{kirby} R. Kirby. \emph{Stable Homeomorphisms and the annulus 
conjecture}. Ann of Math (2) 89, 1969, 575-582.
\bibitem{kulkarni1} R. S. Kulkarni. \emph{Groups with domains of
  discontinuity}.  Math. Ann. 237 (1978), 253-272.
\bibitem{kulkarni} R. S. Kulkarni. \emph{Conformal structures and
M\"obius structures}. Aspects of Mathematics, edited by R.S. Kulkarni
and U. Pinkhall, Max Planck Institut fur Mathematik, Vieweg (1988).
\bibitem{maskit} B. Maskit. \emph{Kleinian Groups}. Springer
Verlag, 1997.
\bibitem{mcmillan} D. R. McMillan Jr., T. L. Thickstun. \emph{Open
three-manifolds and the Poin\-ca\-r\'e Con\-jec\-ture}. Topology 19 (1980),
no. 3, 313-320.
\bibitem{mcmullen} C. T. McMullen. \emph{Renormalization and 3-manifolds
with Fiber over the Circle}. Annals of Mathematics, Studies
142. Princeton University Press, 1996.
\bibitem{milnor} J. Milnor. \emph{Singular points of Complex
Hypersurfaces}. Annales of Mathematics, Studies 61. Princeton
University Press, 1968.
\bibitem{mumford} D. Mumford, C. Series, D. Wright. \emph{Indra's
  Pearls. The vision of Felix Klein}. Cambridge University Press, New York, 2002.
\bibitem{alberto2}Le Dung Trang, J. Seade, A. Verjovsky. 
\emph{Quadrics, Orthogonal Actions and Involutions in 
Complex Proyective Spaces}. To appear (2002).
\bibitem{palais} R. Palais. \emph{Local triviality of the
restriction map for embeddings}. Comment. Math. Helv. 34, 1960, 305-312.
\bibitem{richards} I. Richards. \emph{On the Classification of
Noncompact Surfaces}. Trans. Amer. Math. Soc. 106 (1963), 259-269. 
\bibitem{rolfsen} D. Rolfsen. \emph{Knots and Links}. Publish or
Perish, Inc. 1976.
\bibitem{rushing} B. Rushing. \emph{Topological
Embeddings}. Academic Press, 1973, Vol 52.
\bibitem{alberto1} J. Seade, A. Verjovsky. \emph{Higher 
dimensional complex Kleinian Groups}. Math Ann 322 (2002), No. 2, 279-300.
\bibitem{spivak1} M. Spivak. \emph{A Comprehensive Introduction
to Differential Geometry}. Publish or Perish, Inc. 1970.
\bibitem{tukia}P. Tukia. \emph{On isomorphisms of geometrically finite 
Mobius groups}. Publ. Math. I.H.E.S. Vol. 61 (1985), 171-214.
\bibitem{thurston1}W. P. Thurston. \emph{The geometry and topology of
3-manifolds}. Notes. Princeton University 1976-1979.
\bibitem{thurston} W. P. Thurston. \emph{Three-Dimensional Geometry and
Topology}, volume 1. Princeton Mathematical Series 35, Princeton
University Press, 1997.
\bibitem{verjovsky} A. Verjovsky. \emph{Sistemas de Anosov}. 
Monograf\'{\i}as del IMCA, XII-ELAM. 1999.
\bibitem{winkelnkempe}H. E. Winkelnkemper. \emph{Manifolds as open books}.
Bul. Amer. Math. Soc. Vol. 79 (1973), 45-51.
\bibitem{zeeman1} E. C. Zeeman. \emph{Twisting Spun
Knots}. Trans. Amer. Math. Soc. 115 (1965), 471-495.

\end{thebibliography}
\end{document}